\documentclass[12pt,leqno]{article}
\usepackage[mathscr]{eucal}
\usepackage{amsmath,amssymb,latexsym,theorem,bbm,amsfonts, dsfont}
\usepackage{verbatim}
\usepackage{color}
\usepackage{appendix}
\usepackage{array,epsfig,theorem,bbm,graphics,booktabs,url}

\newcommand{\SC}{\scriptstyle}

\newcommand{\CC}{\mathbb{C}}

\newcommand{\NN}{\mathbb{N}}

\newcommand{\RR}{\mathbb{R}}

\newcommand{\ZZ}{\mathbb{Z}}

\newcommand{\bA}{{\boldsymbol{A}}}

\newcommand{\bB}{{\boldsymbol{B}}}

\newcommand{\bC}{{\boldsymbol{C}}}

\newcommand{\bD}{{\boldsymbol{D}}}

\newcommand{\be}{{\boldsymbol{e}}}

\newcommand{\bI}{{\boldsymbol{I}}}

\newcommand{\bm}{{\boldsymbol{m}}}
\newcommand{\bM}{{\boldsymbol{M}}}

\newcommand{\bP}{{\boldsymbol{P}}}

\newcommand{\bQ}{{\boldsymbol{Q}}}

\newcommand{\bS}{{\boldsymbol{S}}}

\newcommand{\bu}{{\boldsymbol{u}}}

\newcommand{\bv}{{\boldsymbol{v}}}
\newcommand{\bV}{{\boldsymbol{V}}}

\newcommand{\bx}{{\boldsymbol{x}}}
\newcommand{\bX}{{\boldsymbol{X}}}

\newcommand{\bY}{{\boldsymbol{Y}}}

\newcommand{\bZ}{{\boldsymbol{Z}}}
\newcommand{\bU}{{\boldsymbol{U}}}

\newcommand{\bxi}{{\boldsymbol{\xi}}}

\newcommand{\btheta}{{\boldsymbol{\theta}}}
\newcommand{\bvare}{{\boldsymbol{\vare}}}

\newcommand{\bzero}{{\boldsymbol{0}}}
\newcommand{\bone}{{\boldsymbol{1}}}

\newcommand{\cD}{{\mathcal D}}

\newcommand{\cG}{{\mathcal G}}
\newcommand{\cF}{{\mathcal F}}

\newcommand{\cN}{{\mathcal N}}

\newcommand{\cX}{{\mathcal X}}

\newcommand{\bcX}{\boldsymbol{\cX}}

\newcommand{\dd}{\mathrm{d}}

\newcommand{\ff}{\mathrm{f}}
\newcommand{\ii}{\mathrm{i}}

\newcommand{\EE}{\operatorname{\mathbb{E}}}

\newcommand{\PP}{\operatorname{\mathbb{P}}}

\newcommand{\var}{\operatorname{Var}}
\newcommand{\cov}{\operatorname{Cov}}

\newcommand{\vare}{\varepsilon}

\renewcommand{\mid}{\,|\,}

\renewcommand{\leq}{\leqslant}
\renewcommand{\geq}{\geqslant}

\newcommand{\stoch}{\stackrel{\PP}{\longrightarrow}}
\newcommand{\distr}{\stackrel{\cD}{\longrightarrow}}
\newcommand{\distrf}{\stackrel{\cD_\ff}{\longrightarrow}}
\newcommand{\distre}{\stackrel{\cD}{=}}

\newcommand{\as}{\stackrel{{\mathrm{a.s.}}}{\longrightarrow}}

\newcommand{\bbone}{\mathbbm{1}}
\newcommand{\ns}{{\lfloor ns\rfloor}}
\newcommand{\nt}{{\lfloor nt\rfloor}}

\newcommand{\proofend}{\hfill\mbox{$\Box$}}

\numberwithin{equation}{section}

\theoremstyle{change} \theorembodyfont{\em}
\newtheorem{Lem}{Lemma.}[section]
\newtheorem{Thm}[Lem]{Theorem.}
\newtheorem{Pro}[Lem]{Proposition.}

\theorembodyfont{\rm}
\newtheorem{Rem}[Lem]{Remark.}

\begin{document}

\begin{center}
 {\bfseries\Large On aggregation of multitype Galton--Watson \\[1mm] branching processes with immigration}\\[1mm]
 {\bfseries\Large }

\vskip0.5cm

 {\sc\large
  M\'aty\'as $\text{Barczy}^{*,\diamond}$,
  Fanni $\text{K. Ned\'enyi}^{**}$,
  Gyula $\text{Pap}^{**}$}
\end{center}

\vskip0.2cm

\noindent
 * MTA-SZTE Analysis and Stochastics Research Group,
   Bolyai Institute, University of Szeged,
   Aradi v\'ertan\'uk tere 1, H--6720 Szeged, Hungary.

\noindent
 ** Bolyai Institute, University of Szeged,
     Aradi v\'ertan\'uk tere 1, H--6720 Szeged, Hungary.

\noindent e--mails: barczy@math.u-szeged.hu (M. Barczy),
                    nfanni@math.u-szeged.hu (F. K. Ned\'enyi),
                    papgy@math.u-szeged.hu (G. Pap).

\noindent $\diamond$ Corresponding author.

\vskip0.2cm


\renewcommand{\thefootnote}{}
\vspace*{0.2cm}
\footnote{M\'aty\'as Barczy is supported by the J\'anos Bolyai Research Scholarship of the Hungarian Academy of
 Sciences.}

\vspace*{-10mm}

\begin{abstract}
Limit behaviour of temporal and contemporaneous aggregations of independent copies of a stationary multitype
 Galton--Watson branching process with immigration is studied in the so-called iterated and simultaneous cases,
 respectively.
In both cases, the limit process is a zero mean Brownian motion with the same covariance function under third order
 moment conditions on the branching and immigration distributions.
We specialize our results for generalized integer-valued autoregressive processes and single-type Galton--Watson
 processes with immigration as well.
\end{abstract}

\section{Introduction}

The field of temporal and contemporaneous aggregations of independent stationary stochastic processes is an
 important and very active research area in the empirical and theoretical statistics and in other areas as well.
The scheme of contemporaneous (also called cross-sectional) aggregation of random-coefficient autoregressive
 processes of order 1 was firstly proposed by Robinson \cite{Rob} and Granger \cite{Gra} in order to obtain the
 long memory phenomena in aggregated time series.
For surveys on papers dealing with the aggregation of different kinds of stochastic processes, see, e.g.,
 Pilipauskait{\.e} and Surgailis \cite{PilSur}, Jirak \cite[page 512]{Jir} or the arXiv version of Barczy et al.\
 \cite{BarNedPap}.

In this paper we study the limit behaviour of temporal (time) and contemporaneous (space) aggregations of
 independent copies of a strictly stationary multitype Galton--Watson branching process with immigration
 in the so-called iterated and simultaneous cases, respectively.
According to our knowledge, the aggregation of general multitype Galton--Watson branching processes with
 immigration has not been considered in the literature so far.
To motivate the fact that the aggregation of branching processes could be an important topic, now we present an
 interesting and relevant example, where the phenomena of aggregation of this kind of processes may come into play.
A usual Integer-valued AutoRegressive (INAR) process of order 1, \ $(X_k)_{k\geq 0}$, \ can be used to model
 migration, which is quite a big issue nowadays all over the world.
More precisely, given a camp, for all \ $k \geq 0$, \ the random variable \ $X_k$ \ can be interpreted as the
 number of migrants to be present in the camp at time \ $k$, \ and every migrant will stay in the camp with
 probability \ $\alpha \in (0, 1)$ \ indepedently of each other (i.e., with probability \ $1 - \alpha$ \ each
 migrant leaves the camp) and at any time \ $k \geq 1$ \ new migrants may come to the camp.
Given several camps in a country, we may suppose that the corresponding INAR processes of order 1 share the same
 parameter \ $\alpha$ \ and they are independent.
So, the temporal and contemporaneous aggregations of these INAR processes of order 1 is the total usage of the camps
 in terms of the number of migrants in the given country in a given time period, and this quantity may be worth
 studying.

The present paper is organized as follows.
In Section \ref{BRANCHING} we formulate our main results, namely the iterated and simultaneous limit behaviour of
 time- and space-aggregated independent stationary $p$-type Galton--Watson branching processes with immigration is
 described (where \ $p \geq 1$), \ see Theorems \ref{double_aggregation} and \ref{simulataneous_aggregation}.
The limit distributions in these limit theorems coincide, namely, it is a \ $p$-dimensional zero mean Brownian
 motion with a covariance function depending on the expectations and covariances of the offspring and immigration
 distributions.
In the course of the proofs of our results, in Lemma 2.3, we prove that for a subcritical, positively regular
 multitype Galton--Watson branching process with nontrivial immigration, its unique stationary distribution admits
 finite \ $\alpha^\mathrm{th}$ \ moments provided that the branching and immigration distributions have finite
 \ $\alpha^\mathrm{th}$ \ moments, where \ $\alpha \in \{1, 2, 3\}$.
\ In case of \ $\alpha \in \{1, 2\}$, \ Quine \cite{Quine} contains this result, however in case of \ $\alpha = 3$,
 \ we have not found any precise proof in the literature for it, it is something like a folklore, so we decided to
 write down a detailed proof.
As a by-product, we obtain an explicit formula for the third moment in question.
Section \ref{SPECIAL} is devoted to the special case of generalized INAR processes, especially to single-type
 Galton--Watson branching processes with immigration.
All of the proofs can be found in Section \ref{Proofs}.

\section{Aggregation of multitype Galton--Watson branching processes with immigration}
\label{BRANCHING}

Let \ $\ZZ_+$, \ $\NN$, \ $\RR$, \ $\RR_+$, \ and \ $\CC$ \ denote the set of non-negative integers, positive
 integers, real numbers, non-negative real numbers, and complex numbers, respectively.
For all \ $d \in \NN$, \ the \ $d \times d$ \ identity matrix is denoted by \ $\bI_d$.
\ The standard basis in \ $\RR^d$ \ is denoted by \ $\{\be_1, \ldots, \be_d\}$.
\ For \ $\bv \in \RR^d$, \ the Euclidean norm is denoted by \ $\|\bv\|$, \ and for \ $\bA \in \RR^{d\times d}$,
 \ the induced matrix norm is denoted by \ $\|\bA\|$ \ as well (with a little abuse of notation).
All the random variables will be defined on a probability space \ $(\Omega, \cF, \PP)$.

Let \ $(\bX_k = [X_{k,1}, \dots, X_{k,p}]^\top)_{k\in \ZZ_+}$ \ be a \ $p$-type Galton--Watson branching process
 with immigration.
For each \ $k, \ell \in \ZZ_+$ \ and \ $i, j \in \{1, \ldots, p\}$, \ the number of \ $j$-type individuals in the
 \ $k^\mathrm{th}$ \ generation will be denoted by \ $X_{k,j}$, \ the number of \ $j$-type offsprings produced by
 the \ $\ell^\mathrm{th}$ \ individual belonging to type \ $i$ \ of the \ $(k-1)^\mathrm{th}$ \ generation will be
 denoted by \ $\xi^{(i, j)}_{k,\ell}$, \ and the number of immigrants of type \ $i$ \ in the \ $k^\mathrm{th}$
 \ generation will be denoted by \ $\vare^{(i)}_k$.
\ Then we have
 \begin{equation}\label{GWI}
  \bX_k = \sum_{\ell=1}^{X_{k-1,1}}
           \begin{bmatrix}
            \xi^{(1,1)}_{k,\ell} \\
            \vdots \\
            \xi^{(1,p)}_{k,\ell}
           \end{bmatrix}
          + \dots
          + \sum_{\ell=1}^{X_{k-1,p}}
             \begin{bmatrix}
              \xi^{(p,1)}_{k,\ell} \\
              \vdots \\
              \xi^{(p,p)}_{k,\ell}
            \end{bmatrix}
          + \begin{bmatrix}
             \vare^{(1)}_{k} \\
             \vdots \\
             \vare^{(p)}_{k}
            \end{bmatrix} \\
        =: \sum_{i=1}^p \sum_{\ell=1}^{X_{k-1,i}} \bxi^{(i)}_{k,\ell} + \bvare_k
 \end{equation}
 for every \ $k \in \NN$, \ where we define \ $\sum_{\ell=1}^0 := \bzero$.
\ Here \ $\bigl\{\bX_0, \bxi^{(i)}_{k,\ell}, \bvare_k : k, \ell \in \NN, i \in \{1,\ldots,p\}\bigr\}$ \ are
 supposed to be independent \ $\ZZ_+^p$-valued random vectors.
Note that we do not assume independence among the components of these vectors.
Moreover, for all \ $i \in \{1, \ldots, p\}$, \ $\{\bxi^{(i)}, \bxi^{(i)}_{k,\ell} : k, \ell \in \NN\}$ \ and
 \ $\{\bvare, \bvare_k : k \in \NN\}$ \ are supposed to consist of identically distributed random vectors,
 respectively.

Let us introduce the notations \ $\bm_\bvare := \EE(\bvare) \in \RR_+^p$,
 \ $\bM_\bxi := \EE\bigl(\bigl[\bxi^{(1)}, \dots, \bxi^{(p)}\bigr]\bigr) \in \RR_+^{p\times p}$ \ and
 \[
   \bv_{(i,j)}
   := \bigl[\cov(\xi^{(1,i)}, \xi^{(1,j)}), \dots, \cov(\xi^{(p,i)}, \xi^{(p,j)}),
           \cov(\vare^{(i)}, \vare^{(j)})\bigr]^\top \in \RR^{(p+1)\times1}
 \]
 for \ $i, j \in \{1, \dots, p\}$, \ provided that the expectations and covariances in question are finite.
Let \ $\varrho(\bM_\bxi)$ \ denote the spectral radius of \ $\bM_\bxi $, \ i.e., the maximum of the modulus of the
 eigenvalues of \ $\bM_\bxi $.
\ The process \ $(\bX_k)_{k\in\ZZ_+}$ \ is called subcritical, critical or supercritical if \ $\varrho(\bM_\bxi)$
 \ is smaller than \ $1$, \ equal to \ $1$ \ or larger than \ $1$, \ respectively.
The matrix \ $\bM_\bxi$ \ is called primitive if there is a positive integer \ $n \in \NN$ \ such that all the
 entries of \ $\bM_\bxi^n$ \ are positive.
The process \ $(\bX_k)_{k\in\ZZ_+}$ \ is called positively regular if \ $\bM_\bxi$ \ is primitive.
In what follows, we suppose that
 \begin{equation}\label{assumption}
  \begin{aligned}
   &\EE(\bxi^{(i)}) \in \RR_+^p , \quad i \in \{1, \ldots, p\}, \qquad
    \bm_{\bvare} \in \RR_+^p \setminus \{\bzero\} , \\
   &\hspace*{18mm}\rho(\bM_\bxi ) < 1 , \qquad \text{$\bM_\bxi$ \ is primitive.}
  \end{aligned}
 \end{equation}
For further application, we define the matrix
 \begin{align}\label{V_matrix}
  \bV := (V_{i,j})_{i,j=1}^p
      := \left(\bv_{(i,j)}^\top
               \begin{bmatrix}
                (\bI_p - \bM_\bxi)^{-1} \bm_\bvare \\
                1
               \end{bmatrix}\right)_{i,j=1}^p \in \RR^{p \times p} ,
 \end{align}
 provided that the covariances in question are finite.

\begin{Rem}\label{Minv}
Note that the matrix \ $(\bI_p-\bM_\bxi)^{-1}$, \ which appears in \eqref{V_matrix} and throughout the paper,
 exists.
Indeed, \ $\lambda \in \CC$ \ is an eigenvalue of \ $\bI_p - \bM_\bxi$ \ if and only if \ $1 - \lambda$ \ is that
 of \ $\bM_\bxi $.
\ Therefore, since \ $\rho(\bM_\bxi) < 1$, \  all eigenvalues of \ $\bI_p - \bM_\bxi $ \ are non-zero.
This means that \ $\det(\bI_p - \bM_\bxi) \neq 0$, \ so \ $(\bI_p - \bM_\bxi)^{-1}$ \ does exist.
One could also refer to Corollary 5.6.16 and Lemma 5.6.10 in Horn and Johnson \cite{HornJohnson}.
\proofend
\end{Rem}

\begin{Rem}
Note that \ $\bV$ \ is symmetric and positive semidefinite, since \ $\bv_{(i,j)} = \bv_{(j,i)}$,
 \ $i, j \in \{1,\ldots,p\}$, \ and for all \ $\bx\in\RR^p$,
 \begin{align*}
  \bx^\top \bV \bx
  = \sum_{i=1}^p \sum_{j=1}^p V_{i,j} x_i x_j
  = \left(\sum_{i=1}^p \sum_{j=1}^p x_i x_j \bv_{(i,j)}^\top \right)
    \begin{bmatrix}
     (\bI_p - \bM_\bxi)^{-1} \bm_\bvare \\
     1
    \end{bmatrix} ,
 \end{align*}
 where
 \begin{align*}
  \sum_{i=1}^p \sum_{j=1}^p x_i x_j \bv_{(i,j)}^\top
  = \bigl[\bx^\top \cov(\bxi^{(1)}, \bxi^{(1)}) \bx , \ldots, \bx^\top \cov(\bxi^{(p)}, \bxi^{(p)}) \bx ,
          \bx^\top \cov(\bvare, \bvare) \bx \bigr] .
 \end{align*}
Here \ $\bx^\top \cov(\bxi^{(i)}, \bxi^{(i)}) \bx \geq 0$, \ $i \in \{1, \ldots, p\}$,
 \ $\bx^\top \cov(\bvare, \bvare) \bx \geq 0$, \ and \ $(\bI_p - \bM_\bxi)^{-1} \bm_\bvare \in \RR_+^p$ \ due to
 the fact that \ $(\bI_p - \bM_\bxi)^{-1} \bm_\bvare$ \ is nothing else but the expectation vector of the unique
 stationary distribution of \ $(\bX_k)_{k\in\ZZ_+}$, \ see the discussion below and formula \eqref{moment_1_pi}.
\proofend
\end{Rem}

Under \eqref{assumption}, by the Theorem in Quine \cite{Quine}, there is a unique stationary distribution \ $\pi$
 \ for \ $(\bX_k)_{k\in\ZZ_+}$.
\ Indeed, under \eqref{assumption}, \ $\bM_\bxi$ \ is irreducible following from the primitivity of \ $\bM_\bxi$,
 \ see Definition 8.5.0 and Theorem 8.5.2 in Horn and Johnson \cite{HornJohnson}.
For the definition of irreducibility, see Horn and Johnson \cite[Definitions 6.2.21 and 6.2.22]{HornJohnson}.
Further, \ $\bM_\bxi$ \ is aperiodic, since this is equivalent to the primitivity of \ $\bM_\bxi$, \ see Kesten and Stigum
 \cite[page 314]{KesSti3} and Kesten and Stigum \cite[Section 3]{KesSti2}.
For the definition of aperiodicity (also called acyclicity), see, e.g., the Introduction of Danka and Pap
 \cite{DanPap}.
Finally, since \ $\bm_\bvare \in \RR_+^p \setminus \{\bzero\}$, \ the probability generator function of
 \ $\bvare$ \ at \ $\bzero$ \ is less than \ $1$, \ and
 \[
   \EE\Biggl(\log\Biggl(\sum_{i=1}^p \vare^{(i)}\Biggr) \bbone_{\{\bvare\ne\bzero\}}\Biggr)
   \leq \EE\Biggl(\sum_{i=1}^p \vare^{(i)} \bbone_{\{\bvare\ne\bzero\}}\Biggr)
   \leq \EE\left(\sum_{i=1}^p \vare^{(i)}\right)
   = \sum_{i=1}^p \EE(\vare^{(i)}) < \infty ,
 \]
 so one can apply the Theorem in Quine \cite{Quine}.

For each \ $\alpha \in \NN$, \ we say that the \ $\alpha^\mathrm{th}$ \ moment of a random vector is finite if all
 of its mixed moments of order \ $\alpha$ \ are finite.

\begin{Lem}\label{stationary_moments}
Let us assume \eqref{assumption}.
For each \ $\alpha \in \{1, 2, 3\}$, \ the unique stationary distribution \ $\pi$ \ has a finite
 \ $\alpha^\mathrm{th}$ \ moment, provided that the \ $\alpha^\mathrm{th}$ \ moments of \ $\bxi^{(i)}$,
 \ $i \in \{1, \ldots, p\}$, \ and \ $\bvare$ \ are finite.
\end{Lem}

In what follows, we suppose \eqref{assumption} and that the distribution of \ $\bX_0$ \ is the unique
 stationary distribution \ $\pi$, \ hence the Markov chain \ $(\bX_k)_{k\in\ZZ_+}$ \ is strictly stationary.
Recall that, by (2.1) in Quine and Durham \cite{QuineDurham}, for any measurable function
 \ $f : \RR^p \to \RR$ \ satisfying \ $\EE(|f(\bX_0)|) < \infty$, \ we have
 \begin{align}\label{help_ergodic}
  \frac{1}{n} \sum_{k=1}^n f(\bX_k) \as \EE( f(\bX_0)) \qquad \text{as \ $n\to\infty$.}
 \end{align}
First we consider a simple aggregation procedure.
For each \ $N \in \NN$, \ consider the stochastic process \ $\bS^{(N)} = (\bS^{(N)}_k)_{k\in\ZZ_+}$ \ given by
 \[
   \bS^{(N)}_k := \sum_{j=1}^N (\bX^{(j)}_k - \EE(\bX^{(j)}_k)) , \qquad k \in \ZZ_+ ,
 \]
 where \ $\bX^{(j)} = (\bX^{(j)}_k)_{k\in\ZZ_+}$, \ $j \in \NN$, \ is a sequence of independent copies of the
 strictly stationary $p$-type Galton--Watson process \ $(\bX_k)_{k\in\ZZ_+}$ \ with immigration.
Here we point out that we consider so-called idiosyncratic immigrations, i.e., the immigrations belonging to
 \ $\bX^{(j)}$, \ $j\in\NN$, \ are independent.

We will use \ $\distrf$ \ or \ $\cD_\ff\text{-}\hspace*{-1mm}\lim$ \ for weak convergence of finite dimensional
 distributions, and \ $\distr$ \ for weak convergence in \ $D(\RR_+, \RR^p)$ \ of stochastic processes with
 c\`adl\`ag sample paths, where \ $D(\RR_+, \RR^p)$ \ denotes the space of $\RR^p$-valued c\`adl\`ag functions
 defined on \ $\RR_+$.

\begin{Pro}\label{simple_aggregation}
If all entries of the vectors \  $\bxi^{(i)}$, \ $i \in \{1, \ldots, p\}$, \ and \ $\bvare$ \ have finite second
 moments, then
 \[
   N^{-\frac{1}{2}} \bS^{(N)} \distrf \bcX \qquad \text{as \ $N \to \infty$,}
 \]
 where \ $\bcX = (\bcX_k)_{k\in\ZZ_+}$ \ is a stationary \ $p$-dimensional zero mean Gaussian process with
 covariances
 \begin{align}\label{cov}
  \EE(\bcX_0 \bcX_k^\top)
  = \cov(\bX_0, \bX_k)
  = \var(\bX_0) (\bM_\bxi^\top)^k, \qquad k \in \ZZ_+ ,
 \end{align}
 where
 \begin{equation}\label{var0}
  \var(\bX_0) = \sum_{k=0}^\infty \bM_\bxi^k \bV (\bM_\bxi^{\top})^k .
 \end{equation}
\end{Pro}

We note that using formula \eqref{moment_2_pi} presented later on, one could give an explicit formula for
 \ $\var(\bX_0)$ \ (not containing an infinite series).

\begin{Pro}\label{simple_aggregation2}
If all entries of the vectors \  $\bxi^{(i)}$, \ $i \in \{1, \ldots, p\}$, \ and \ $\bvare$ \ have finite third
 moments, then
 \[
   \biggl(n^{-\frac{1}{2}} \sum_{k=1}^\nt \bS^{(1)}_k\biggr)_{t\in\RR_+}
   = \biggl(n^{-\frac{1}{2}} \sum_{k=1}^\nt (\bX_k^{(1)} - \EE(\bX_k^{(1)}))\biggr)_{t\in\RR_+}
   \distr (\bI_p - \bM_\bxi)^{-1} \bB
   \qquad \text{as \ $n \to \infty$,}
 \]
 where \ $\bB = (\bB_t)_{t\in\RR_+}$ \ is a \ $p$-dimensional zero mean Brownian motion satisfying
 \ $\var(\bB_1) = \bV$.
\end{Pro}

Note that Propositions \ref{simple_aggregation} and \ref{simple_aggregation2} are about the scalings of the
 space-aggregated process \ $\bS^{(N)}$ \ and the time-aggregated process
 \ $\bigl(\sum_{k=1}^\nt \bS^{(1)}_k\bigr)_{t\in\RR_+}$, \ respectively.

For each \ $N, n \in \NN$, \ consider the stochastic process \ $\bS^{(N,n)} = (\bS^{(N,n)}_t)_{t\in\RR_+}$ \ given
 by
 \[
   \bS^{(N,n)}_t := \sum_{j=1}^N \sum_{k=1}^\nt (\bX^{(j)}_k - \EE(\bX^{(j)}_k)) , \qquad t \in \RR_+ .
 \]


\begin{Thm}\label{double_aggregation}
If all entries of the vectors \ $\bxi^{(i)}$, \ $i \in \{1,\ldots,p\}$, \ and \ $\bvare$ \ have finite second
 moments, then
 \begin{align}\label{help1_double_aggregation}
  \cD_\ff\text{-}\hspace*{-1mm}\lim_{n\to\infty} \,
  \cD_\ff\text{-}\hspace*{-1mm}\lim_{N\to\infty} \,
   (nN)^{-\frac{1}{2}} \bS^{(N,n)}
  = (\bI_p - \bM_\bxi)^{-1} \bB ,
 \end{align}
 where \ $\bB = (\bB_t)_{t\in\RR_+}$ \ is a \ $p$-dimensional zero mean Brownian motion satisfying
 \ $\var(\bB_1) = \bV$.
\ If all entries of the vectors \  $\bxi^{(i)}$, \ $i \in \{1, \ldots, p\}$, \ and \ $\bvare$ \ have finite
 third moments, then
 \begin{align}\label{help2_double_aggregation}
  \cD_\ff\text{-}\hspace*{-1mm}\lim_{N\to\infty} \,
  \cD_\ff\text{-}\hspace*{-1mm}\lim_{n\to\infty} \,
   (nN)^{-\frac{1}{2}} \bS^{(N,n)}
  = (\bI_p - \bM_\bxi)^{-1} \bB ,
 \end{align}
 where \ $\bB = (\bB_t)_{t\in\RR_+}$ \ is a \ $p$-dimensional zero mean Brownian motion satisfying
 \ $\var(\bB_1) = \bV$.
\end{Thm}

\begin{Thm}\label{simulataneous_aggregation}
If all entries of the vectors \ $\bxi^{(i)}$, \ $i \in \{1,\ldots,p\}$, \ and \ $\bvare$ \ have finite third
 moments, then
 \begin{align}\label{help3_simulataneous_aggregation}
  (nN)^{-\frac{1}{2}} \bS^{(N,n)} \distr (\bI_p - \bM_\bxi)^{-1} \bB ,
 \end{align}
 if both \ $n$ \ and \ $N$ \ converge to infinity (at any rate), where \ $\bB = (\bB_t)_{t\in\RR_+}$ \ is a
 \ $p$-dimensional zero mean Brownian motion satisfying \ $\var(\bB_1) = \bV$.
\end{Thm}

A key ingredient of the proofs is the fact that \ $(\bX_k - \EE(\bX_k))_{k\in\ZZ_+}$ \ can be rewritten as a stable
 first order vector autoregressive process with coefficient matrix \ $\bM_\bxi$ \ and with heteroscedastic innovations,
 see \eqref{help5}.

\section{A special case: aggregation of GINAR processes}\label{SPECIAL}

We devote this section to the analysis of aggregation of Generalized Integer-Valued Autoregressive processes of
 order \ $p \in \NN$ \ (GINAR($p$) processes), which are special cases of \ $p$-type Galton--Watson branching
 processes with immigration introduced in \eqref{GWI}.
For historical fidelity, we note that it was Latour \cite{Latour} who introduced GINAR($p$) processes as
 generalizations of INAR($p$) processes.
This class of processes became popular in modelling integer-valued time series data such as the daily number of
 claims at an insurance company.
In fact, a GINAR(1) process is a (general) single type Galton--Watson branching processes with immigration.

Let \ $(Z_k)_{k\geq -p+1}$ \ be a GINAR($p$) process.
Namely, for each \ $k, \ell \in \ZZ_+$ \ and \ $i \in \{1,\dots,p\}$, \ the number of individuals in the
 \ $k^\mathrm{th}$ \ generation will be denoted by \ $Z_k$, \ the number of offsprings produced by the
 \ $\ell^\mathrm{th}$ \ individual belonging to the \ $(k-i)^\mathrm{th}$ \ generation will be denoted by
 \ $\xi^{(i,1)}_{k,\ell}$, \ and the number of immigrants in the \ $k^\mathrm{th}$ \ generation will be denoted by
 \ $\vare^{(1)}_k$.
\ Here the \ $1$-s in the supercripts of \ $\xi^{(i,1)}_{k,\ell}$ \ and \ $\vare^{(1)}_k$ \ are displayed in order
 to have a better comparison with \eqref{GWI}.
Then we have
 \begin{equation*}
  Z_k
  = \sum_{\ell=1}^{Z_{k-1}} \xi^{(1,1)}_{k,\ell} + \dots + \sum_{\ell=1}^{Z_{k-p}} \xi^{(p,1)}_{k,\ell}
    + \vare_k^{(1)} , \qquad k \in \NN .
 \end{equation*}
Here
 \ $\bigl\{Z_0, Z_{-1},\ldots, Z_{-p+1}, \xi^{(i,1)}_{k,\ell}, \vare_k^{(1)}
           : k, \ell \in \NN , i \in \{1, \ldots, p\}\bigr\}$
 \ are supposed to be independent nonnegative integer-valued random variables.
Moreover, for all \ $i \in \{1, \ldots, p\}$, \ $\{\xi^{(i,1)}, \xi^{(i,1)}_{k,\ell} : k, \ell \in \NN \}$ \ and
 \ $\{\vare^{(1)}, \vare_k^{(1)} : k \in \NN \}$ \ are supposed to consist of identically distributed random
 variables, respectively.

A GINAR($p$) process can be embedded in a \ $p$-type Galton--Watson branching process with immigration
 \ $(\bX_k = [Z_k, \dots, Z_{k-p+1}]^\top)_{k\in\ZZ_+}$ \ with the corresponding \ $p$-dimensional random vectors
 \[
   \bxi^{(1)}_{k,\ell} = \begin{bmatrix}
                          \xi^{(1,1)}_{k,\ell} \\
                          1 \\
                          0 \\
                          \vdots \\
                          0
                         \end{bmatrix} ,
   \quad \cdots, \quad
   \bxi^{(p-1)}_{k,\ell} = \begin{bmatrix}
                            \xi^{(p-1,1)}_{k,\ell} \\
                            0 \\
                            \vdots \\
                            0 \\
                            1
                            \end{bmatrix} ,
   \quad
   \bxi^{(p)}_{k,\ell} = \begin{bmatrix}
                          \xi^{(p,1)}_{k,\ell} \\
                          0 \\
                          0 \\
                          \vdots \\
                          0
                         \end{bmatrix} ,
   \quad
   \bvare_k = \begin{bmatrix}
               \vare^{(1)}_k \\
               0 \\
               0 \\
               \vdots \\
               0
              \end{bmatrix}
 \]
 for any \ $k, \ell \in \NN$.

In what follows, we reformulate the classification of GINAR($p$) processes in terms of the expectations of
 the offspring distributions.

\begin{Rem}
In case of a GINAR($p$) process, one can show that \ $\varphi$, \ the characteristic polynomial of the matrix
 \ $\bM_\bxi $, \ has the form
 \[
   \varphi(\lambda)
   := \det(\lambda \bI_p - \bM_\bxi)
    = \lambda^p - \EE(\xi^{(1,1)}) \lambda^{p-1} - \cdots - \EE(\xi^{(p-1,1)}) \lambda - \EE(\xi^{(p,1)}) ,
   \qquad \lambda \in \CC .
 \]
Recall that \ $\varrho(\bM_\bxi)$ \ denotes the spectral radius of \ $\bM_\bxi $, \ i.e., the maximum of the
 modulus of the eigenvalues of \ $\bM_\bxi $.
\ If \ $\EE(\xi^{(p,1)}) > 0$, \ then, by the proof of Proposition 2.2 in Barczy et al. \cite{BarIspPap2}, the
 characteristic polynomial \ $\varphi$ \ has just one positive root, \ $\varrho(\bM_\bxi) > 0$, \ the nonnegative
 matrix \ $\bM_\bxi$ \ is irreducible, \ $\varrho(\bM_\bxi)$ \ is an eigenvalue of \ $\bM_\bxi$, \ and
 \ $\sum_{i=1}^p \EE(\xi^{(i,1)}) \varrho(\bM_\bxi)^{-i} = 1$.
\ Further,
 \[
   \varrho(\bM_\bxi ) \,
   \begin{cases}
    < & \\
    = & \\
    >
   \end{cases}
   1
   \qquad \Longleftrightarrow  \qquad
   \sum_{i=1}^p \EE(\xi^{(i,1)}) \,
   \begin{cases}
    < & \\
    = & \\
    >
   \end{cases}
   1 . \\[-4mm]
 \]
\proofend
\end{Rem}

Next, we specialize the matrix \ $\bV$, \ defined in \eqref{V_matrix}, in case of a subcritical GINAR($p$)
 process.

\begin{Rem}
In case of a GINAR($p$) process, the vectors
 \[
   \bv_{(i,j)}
   = \bigl[\cov(\xi^{(1,i)}, \xi^{(1,j)}), \dots, \cov(\xi^{(p,i)}, \xi^{(p,j)}),
           \cov(\vare^{(i)}, \vare^{(j)})\bigr]^\top
   \in \RR^{(p+1)\times1}
 \]
 for \ $i, j \in\{1, \dots, p\}$ \ are all zero vectors except for the case \ $i = j = 1$.
\ Therefore, in case of \ $\varrho(\bM_\bxi) < 1$, \ the matrix \ $\bV$, \ defined in \eqref{V_matrix}, reduces to
 \begin{equation}\label{V_GINAR}
  \bV = \bv_{(1,1)}^\top
        \begin{bmatrix}
         (\bI_p - \bM_\bxi)^{-1} \EE(\vare^{(1)}) \be_1 \\
         1
        \end{bmatrix}
        (\be_1 \be_1^\top) .
 \end{equation}
\proofend
\end{Rem}

Finally, we specialize the limit distribution in Theorems \ref{double_aggregation} and
 \ref{simulataneous_aggregation} in case of a subcritical GINAR($p$) process.

\begin{Rem}
Let us note that in case of \ $p = 1$ \ and \ $\EE(\xi^{(1,1)}) < 1$ \ (yielding that the corresponding
 GINAR($1$) process is subcritical), the limit process in Theorems \ref{double_aggregation} and
 \ref{simulataneous_aggregation} can be written as
 \[
   \frac{1}{1-\EE(\xi^{(1,1)})}
   \sqrt{\frac{\EE(\vare^{(1)}) \var(\xi^{(1,1)}) + (1 - \EE(\xi^{(1,1)})) \var(\vare^{(1)})}
        {1-\EE(\xi^{(1,1)})}}
   W ,
 \]
 where \ $W = (W_t)_{t\in\RR_+}$ \ is a standard 1-dimensional Brownian motion.
Indeed, this holds, since in this special case \ $\bM_\bxi = \EE(\xi^{(1,1)})$ \ yielding that
 \ $(\bI_p - \bM_\bxi)^{-1} = (1 - \EE(\xi^{(1,1)}))^{-1}$, \ and, by \eqref{V_GINAR},
 \begin{align*}
  \bV = \begin{bmatrix}
         \cov(\xi^{(1,1)}, \xi^{(1,1)}) \\ \cov(\vare^{(1)}, \vare^{(1)})
        \end{bmatrix}^\top
        \begin{bmatrix}
         \frac{\EE(\vare^{(1)})} {1 - \EE(\xi^{(1,1)})} \\
         1 \\
        \end{bmatrix}
      = \frac{\var(\xi^{(1,1)})\EE(\vare^{(1)})}{1-\EE(\xi^{(1,1)})} + \var(\vare^{(1)}) .\\[-17mm]
 \end{align*}
\proofend
\end{Rem}

\section{Proofs}
\label{Proofs}

\noindent{\bf Proof of Lemma \ref{stationary_moments}.}
Let \ $(\bZ_k)_{k\in\ZZ_+}$ \ be a \ $p$-type Galton--Watson branching process without immigration, with the same
 offspring distribution as \ $(\bX_k)_{k\in\ZZ_+}$, \ and with  \ $\bZ_0 \distre \bvare$.
\ Then the stationary distribution \ $\pi$ \ of \ $(\bX_k)_{k\in\ZZ_+}$ \ admits the representation
 \[
   \pi \distre \sum_{r=0}^\infty \bZ_r^{(r)} ,
 \]
 where \ $(\bZ_k^{(n)})_{k\in\ZZ_+}$, \ $n \in \ZZ_+$, \ are independent copies of \ $(\bZ_k)_{k\in\ZZ_+}$.
\ This is a consequence of formula (16) for the probability generating function of \ $\pi$ \ in Quine \cite{Quine}.
It is convenient to calculate moments of Kronecker powers of random vectors.
We will use the notation \ $\bA \otimes \bB$ \ for the Kronecker product of the matrices \ $\bA$ \ and \ $\bB$,
 \ and we put \ $\bA^{\otimes2} := \bA \otimes \bA$ \ and \ $\bA^{\otimes3} := \bA \otimes \bA \otimes \bA$.
\ For each \ $\alpha \in \{1, 2, 3\}$, \ by the monotone convergence theorem, we have
 \[
   \int_{\RR^p} \bx^{\otimes\alpha} \, \pi(\dd\bx)
   = \EE\Biggl[\Biggl(\sum_{r=0}^\infty \bZ_r^{(r)}\Biggr)^{\otimes\alpha}\Biggr]
   = \lim_{n\to\infty} \EE\Biggl[\Biggl(\sum_{r=0}^{n-1} \bZ_r^{(r)}\Biggr)^{\otimes\alpha}\Biggr] .
 \]
For each \ $n \in \ZZ_+$, \ we have
 \[
   \sum_{r=0}^{n-1} \bZ_r^{(r)} \distre \bY_n ,
 \]
 where \ $(\bY_k)_{k\in\ZZ_+}$ \ is a Galton--Watson branching process with the same offspring and immigration
 distributions as \ $(\bX_k)_{k\in\ZZ_+}$, \ and with  \ $\bY_0 = \bzero$.
\ This can be checked comparing their probability generating functions taking into account formula (3) in Quine
 \cite{Quine} as well.
Consequently, we conclude
 \begin{equation}\label{reprpi}
  \int_{\RR^p} \bx^{\otimes\alpha} \, \pi(\dd\bx) = \lim_{n\to\infty} \EE\bigl(\bY_n^{\otimes\alpha}\bigr) .
 \end{equation}
For each \ $n \in \NN$, \ using \eqref{GWI}, we obtain
 \begin{align}\label{help12}
  \begin{split}
  \EE(\bY_n \mid \cF_{n-1}^\bY)
  &= \sum_{i=1}^p \sum_{j=1}^{Y_{n-1,i}} \EE(\bxi_{n,j}^{(i)} \mid \cF_{n-1}^\bY)
     + \EE(\bvare_n \mid \cF_{n-1}^\bY)
   = \sum_{i=1}^p Y_{n-1,i} \EE(\bxi^{(i)}) + \EE(\bvare) \\
  &= \sum_{i=1}^p \EE(\bxi^{(i)}) \be_i^\top \bY_{n-1} + \bm_\bvare
   = \bM_\bxi \bY_{n-1} + \bm_\bvare ,
  \end{split}
 \end{align}
 where \ $\cF_{n-1}^\bY  := \sigma(\bY_0, \ldots, \bY_{n-1})$, \ $n \in \NN$, \ and
 \ $Y_{n-1,i} := \be_i^\top \bY_{n-1}$, \ $i \in \{1, \ldots, p\}$.
\ Taking the expectation, we get
 \begin{equation}\label{moment_1_Y}
  \EE(\bY_n) = \bM_\bxi \EE(\bY_{n-1}) + \bm_\bvare , \qquad n \in \NN.
 \end{equation}
Taking into account \ $\bY_0 = \bzero$, \ we obtain
 \[
   \EE(\bY_n)
   = \sum_{k=1}^n \bM_\bxi^{n-k} \bm_\bvare
   = \sum_{\ell=0}^{n-1} \bM_\bxi^\ell \bm_\bvare , \qquad n \in \NN .
 \]
For each \ $n \in \NN$, \ we have  \ $(\bI_p - \bM_\bxi) \sum_{\ell=0}^{n-1} \bM_\bxi^\ell = \bI_p - \bM_\bxi^n$.
\ By the condition \ $\varrho(\bM_\bxi) < 1$, \ the matrix \ $\bI_p - \bM_\bxi$ \ is invertible and
 \ $\sum_{\ell=0}^\infty \bM_\bxi^\ell = (\bI_p - \bM_\bxi)^{-1}$, \ see Corollary 5.6.16 and Lemma 5.6.10 in Horn
 and Johnson \cite{HornJohnson}.
Consequently, by \eqref{reprpi}, the first moment of \ $\pi$ \ is finite, and
 \begin{equation}\label{moment_1_pi}
  \int_{\RR^p} \bx \, \pi(\dd\bx) = (\bI_p - \bM_\bxi)^{-1} \bm_\bvare .
 \end{equation}

Now we suppose that the second moments of \ $\bxi^{(i)}$, \ $i \in \{1, \ldots, p\}$, \ and \ $\bvare$ \ are
 finite.
For each \ $n \in \NN$, \ using again \eqref{GWI}, we obtain
 \begin{align*}
  \EE(\bY_n^{\otimes2} \mid \cF_{n-1}^\bY)
  &= \sum_{i=1}^p \sum_{j=1}^{Y_{n-1,i}} \sum_{i'=1}^p \sum_{j'=1}^{Y_{n-1,i'}}
      \EE(\bxi_{n,j}^{(i)} \otimes \bxi_{n,j'}^{(i')} \mid \cF_{n-1}^\bY) \\
  &\quad
     + \sum_{i=1}^p \sum_{j=1}^{Y_{n-1,i}}
        \EE(\bxi_{n,j}^{(i)} \otimes \bvare_n + \bvare_n \otimes \bxi_{n,j}^{(i)} \mid \cF_{n-1}^\bY)
     + \EE(\bvare_n^{\otimes2} \mid \cF_{n-1}^\bY) \\
  &= \sum_{i=1}^p \sum_{\underset{\SC i'\ne i}{i'=1}}^p
      Y_{n-1,i} Y_{n-1,i'} \EE(\bxi^{(i)}) \otimes \EE(\bxi^{(i')})
     + \sum_{i=1}^p Y_{n-1,i} (Y_{n-1,i} -1) [\EE(\bxi^{(i)})]^{\otimes2} \\
  &\quad
     + \sum_{i=1}^p Y_{n-1,i} \EE[(\bxi^{(i)})^{\otimes2}]
     + \sum_{i=1}^p Y_{n-1,i} \EE(\bxi^{(i)} \otimes \bvare + \bvare \otimes \bxi^{(i)})
     + \EE(\bvare^{\otimes2})
 \end{align*}
 \begin{align*}
  &= \sum_{i=1}^p \sum_{i'=1}^p Y_{n-1,i} Y_{n-1,i'} \EE(\bxi^{(i)}) \otimes \EE(\bxi^{(i')})
     + \sum_{i=1}^p Y_{n-1,i} \bigl\{\EE[(\bxi^{(i)})^{\otimes2}] - [\EE(\bxi^{(i)})]^{\otimes2}\bigr\} \\
  &\quad
     + \sum_{i=1}^p Y_{n-1,i} \bigl\{\EE(\bxi^{(i)}) \otimes \bm_\bvare + \bm_\bvare \otimes \EE(\bxi^{(i)})\bigr\}
     + \EE(\bvare^{\otimes2}) \\
  &= (\bM_\bxi \bY_{n-1})^{\otimes2} + \bA_{2,1} \bY_{n-1} + \EE(\bvare^{\otimes2}) .
 \end{align*}
 with
 \[
   \bA_{2,1}
   := \sum_{i=1}^p
       \bigl\{\EE[(\bxi^{(i)})^{\otimes2}] + \EE(\bxi^{(i)}) \otimes \bm_\bvare
              + \bm_\bvare \otimes \EE(\bxi^{(i)}) - [\EE(\bxi^{(i)})]^{\otimes2}\bigr\}
       \be_i^\top
   \in \RR^{p^2\times p} .
 \]
Indeed, using the mixed-product property \ $(\bA \otimes \bB) (\bC \otimes \bD) = (\bA \bC) \otimes (\bB \bD)$
 \ for matrices of such size that one can form the matrix products \ $\bA \bC$ \ and \ $\bB \bD$, \ we have
 \[
   Y_{n-1,i} Y_{n-1,i'}
   = Y_{n-1,i} \otimes Y_{n-1,i'}
   = (\be_i^\top \bY_{n-1}) \otimes (\be_{i'}^\top \bY_{n-1})
   = (\be_i^\top \otimes \be_{i'}^\top) \bY_{n-1}^{\otimes2} ,
 \]
 hence
 \begin{align*}
  &\sum_{i=1}^p \sum_{i'=1}^p Y_{n-1,i} Y_{n-1,i'} \EE(\bxi^{(i)}) \otimes \EE(\bxi^{(i')})
   = \sum_{i=1}^p \sum_{i'=1}^p
      \bigl[\EE(\bxi^{(i)}) \otimes \EE(\bxi^{(i')})\bigr] (\be_i^\top \otimes \be_{i'}^\top) \bY_{n-1}^{\otimes2} \\
  &= \sum_{i=1}^p \sum_{i'=1}^p
      \bigl[(\EE(\bxi^{(i)}) \be_i^\top) \otimes (\EE(\bxi^{(i')}) \be_{i'}^\top)\bigr] \bY_{n-1}^{\otimes2}
   = \Biggl(\sum_{i=1}^p \EE(\bxi^{(i)}) \be_i^\top\Biggr)^{\otimes2} \bY_{n-1}^{\otimes2} \\
  &= (\bM_\bxi)^{\otimes2} \bY_{n-1}^{\otimes2}
   = (\bM_\bxi \bY_{n-1})^{\otimes2} .
 \end{align*}
Consequently, we obtain
 \[
   \EE(\bY_n^{\otimes2} \mid \cF_{n-1}^\bY)
   = \bM_\bxi^{\otimes2} \bY_{n-1}^{\otimes2} + \bA_{2,1} \bY_{n-1} + \EE(\bvare^{\otimes2}) ,
   \qquad n \in \NN .
 \]
Taking the expectation, we get
 \begin{equation}\label{moment_2_Y}
  \EE(\bY_n^{\otimes2})
  = \bM_\bxi^{\otimes2} \EE(\bY_{n-1}^{\otimes2}) + \bA_{2,1} \EE(\bY_{n-1}) + \EE(\bvare^{\otimes2}) ,
  \qquad n \in \NN .
 \end{equation}
Using also \eqref{moment_1_Y}, we obtain
 \[
   \begin{bmatrix}
    \EE(\bY_n) \\
    \EE(\bY_n^{\otimes2})
   \end{bmatrix}
   = \bA_2
     \begin{bmatrix}
      \EE(\bY_{n-1}) \\
      \EE(\bY_{n-1}^{\otimes2})
     \end{bmatrix}
     + \begin{bmatrix}
        \bm_\bvare \\
        \EE(\bvare^{\otimes2})
       \end{bmatrix} , \qquad n \in \NN ,
 \]
 with
 \[
   \bA_2 := \begin{bmatrix}
             \bM_\bxi & \bzero \\
             \bA_{2,1} & \bM_\bxi^{\otimes2}
            \end{bmatrix}
   \in \RR^{(p+p^2)\times(p+p^2)} .
 \]
Taking into account \ $\bY_0 = \bzero$, \ we obtain
 \[
   \begin{bmatrix}
    \EE(\bY_n) \\
    \EE(\bY_n^{\otimes2})
   \end{bmatrix}
   = \sum_{k=1}^n
      \bA_2^{n-k}
      \begin{bmatrix}
       \bm_\bvare \\
       \EE(\bvare^{\otimes2})
      \end{bmatrix}
   = \sum_{\ell=0}^{n-1}
      \bA_2^\ell
      \begin{bmatrix}
       \bm_\bvare \\
       \EE(\bvare^{\otimes2})
      \end{bmatrix} , \qquad n \in \NN .
 \]
We have \ $\varrho(\bA_2) = \max\{\varrho(\bM_\bxi), \varrho(\bM_\bxi^{\otimes2})\}$, \ where
 \ $\varrho(\bM_\bxi^{\otimes2}) = [\varrho(\bM_\bxi)]^2$.
\ Taking into account \ $\varrho(\bM_\bxi) < 1$, \ we conclude \ $\varrho(\bA_2) = \varrho(\bM_\bxi) < 1$, \ and,
 by \eqref{reprpi}, the second moment of \ $\pi$ \ is finite,
 and
 \begin{equation}\label{moment_2_pi}
  \begin{bmatrix}
   \int_{\RR^p} \bx \, \pi(\dd\bx) \\
   \int_{\RR^p} \bx^{\otimes2} \, \pi(\dd\bx) \\
  \end{bmatrix}
  = (\bI_{p+p^2} - \bA_2)^{-1}
    \begin{bmatrix}
     \bm_\bvare \\
     \EE(\bvare^{\otimes2})
    \end{bmatrix} .
 \end{equation}
Since
 \[
   (\bI_{p+p^2} - \bA_2)^{-1}
   = \begin{bmatrix}
      (\bI_p - \bM_\bxi)^{-1} & \bzero \\
      (\bI_{p^2}-\bM_\bxi^{\otimes2})^{-1} \bA_{2,1} (\bI_p - \bM_\bxi)^{-1}
       & (\bI_{p^2}-\bM_\bxi^{\otimes2})^{-1} \\
     \end{bmatrix},
 \]
 we have
 \begin{align*}
  \int_{\RR^p} \bx^{\otimes2} \, \pi(\dd\bx)
    = (\bI_{p^2}-\bM_\bxi^{\otimes2})^{-1} \bA_{2,1}  (\bI_p - \bM_\bxi)^{-1} \bm_\bvare
      + (\bI_{p^2}-\bM_\bxi^{\otimes2})^{-1}
        \EE(\bvare^{\otimes2}).
 \end{align*}

Now we suppose that the third moments of \ $\bxi^{(i)}$, \ $i \in \{1, \ldots, p\}$, \ and \ $\bvare$ \ are finite.
For each \ $n \in \NN$, \ using again \eqref{GWI}, we obtain
 \[
   \EE(\bY_n^{\otimes3} \mid \cF_{n-1}^\bY)
   = S_{n,1} + S_{n,2} + S_{n,3} + \EE(\bvare_n^{\otimes3} \mid \cF_{n-1}^\bY)
 \]
 with
 \begin{align*}
  S_{n,1}
  &:= \sum_{i=1}^p \sum_{j=1}^{Y_{n-1,i}} \sum_{i'=1}^p \sum_{j'=1}^{Y_{n-1,i'}} \sum_{i''=1}^p
       \sum_{j''=1}^{Y_{n-1,i''}}
       \EE(\bxi_{n,j}^{(i)} \otimes \bxi_{n,j'}^{(i')} \otimes \bxi_{n,j''}^{(i'')} \mid \cF_{n-1}^\bY) , \\
  S_{n,2}
  &:= \sum_{i=1}^p \sum_{j=1}^{Y_{n-1,i}} \sum_{i'=1}^p \sum_{j'=1}^{Y_{n-1,i'}}
       \EE(\bxi_{n,j}^{(i)} \otimes \bxi_{n,j'}^{(i')} \otimes \bvare_n
           + \bxi_{n,j}^{(i)} \otimes \bvare_n \otimes \bxi_{n,j'}^{(i')}
           + \bvare_n \otimes \bxi_{n,j}^{(i)} \otimes \bxi_{n,j'}^{(i')} \mid \cF_{n-1}^\bY) , \\
  S_{n,3}
  &:= \sum_{i=1}^p \sum_{j=1}^{Y_{n-1,i}}
       \EE(\bxi_{n,j}^{(i)} \otimes \bvare_n^{\otimes2}
           + \bvare_n \otimes \bxi_{n,j}^{(i)} \otimes \bvare_n
           + \bvare_n^{\otimes2} \otimes \bxi_{n,j}^{(i)} \mid \cF_{n-1}^\bY) .
 \end{align*}
We have
 \begin{align*}
  &S_{n,1}
   = \sum_{i=1}^p \sum_{\underset{\SC i'\ne i}{i'=1}}^p \sum_{\underset{\SC i''\notin\{i,i'\}}{i''=1}}^p
      Y_{n-1,i} Y_{n-1,i'} Y_{n-1,i''} \EE(\bxi^{(i)}) \otimes \EE(\bxi^{(i')}) \otimes \EE(\bxi^{(i'')}) \\
  &+ \sum_{i=1}^p \sum_{\underset{\SC i'\ne i}{i'=1}}^p
      Y_{n-1,i} (Y_{n-1,i} -1) Y_{n-1,i'} \\
  &\phantom{+ \sum_{i=1}^p \sum_{\underset{\SC i'\ne i}{i'=1}}^p}
      \times
      \bigl\{[\EE(\bxi^{(i)})]^{\otimes2} \otimes \EE(\bxi^{(i')})
             + \EE(\bxi^{(i)}) \otimes \EE(\bxi^{(i')}) \otimes \EE(\bxi^{(i)})
             + \EE(\bxi^{(i')}) \otimes [\EE(\bxi^{(i)})]^{\otimes2}\bigr\} \\
  &+ \sum_{i=1}^p \sum_{\underset{\SC i'\ne i}{i'=1}}^p
      Y_{n-1,i} Y_{n-1,i'}
      \bigl\{\EE[(\bxi^{(i)})^{\otimes2}] \otimes \EE(\bxi^{(i')})
             + \EE(\bxi^{(i)} \otimes \bxi^{(i')} \otimes \bxi^{(i)})
             + \EE(\bxi^{(i')}) \otimes \EE[(\bxi^{(i)})^{\otimes2}]\bigr\}
 \end{align*}
 \begin{align*}
  &+ \sum_{i=1}^p Y_{n-1,i} (Y_{n-1,i} -1) (Y_{n-1,i} -2) [\EE(\bxi^{(i)})]^{\otimes3}
   + \sum_{i=1}^p Y_{n-1,i} \EE[(\bxi^{(i)})^{\otimes3}] \\
  &+ \sum_{i=1}^p
      Y_{n-1,i} (Y_{n-1,i} -1)
      \bigl\{\EE[(\bxi^{(i)})^{\otimes2}] \otimes \EE(\bxi^{(i)})
             + \EE(\bxi_{1,1}^{(i)} \otimes \bxi_{1,2}^{(i)} \otimes \bxi_{1,1}^{(i)})
             + \EE(\bxi^{(i)}) \otimes \EE[(\bxi^{(i)})^{\otimes2}]\bigr\} ,
 \end{align*}
 which can be written in the form
 \begin{align*}
  S_{n,1}
  &= \sum_{i=1}^p \sum_{i'=1}^p \sum_{i''=1}^p
      Y_{n-1,i} Y_{n-1,i'} Y_{n-1,i''} \EE(\bxi^{(i)}) \otimes \EE(\bxi^{(i')}) \otimes \EE(\bxi^{(i'')}) \\
  &\quad
     + \sum_{i=1}^p \sum_{i'=1}^p
        Y_{n-1,i} Y_{n-1,i'}
        \bigl\{\EE[(\bxi^{(i)})^{\otimes2}] \otimes \EE(\bxi^{(i')})
               + \EE(\bxi^{(i)} \otimes \bxi^{(i')} \otimes \bxi^{(i)}) \\
  &\phantom{\quad + \sum_{i=1}^p \sum_{i'=1}^p Y_{n-1,i} Y_{n-1,i'} \bigl\{}
               + \EE(\bxi^{(i')}) \otimes \EE[(\bxi^{(i)})^{\otimes2}]
               - [\EE(\bxi^{(i)})]^{\otimes2} \otimes \EE(\bxi^{(i')}) \\
  &\phantom{\quad + \sum_{i=1}^p \sum_{i'=1}^p Y_{n-1,i} Y_{n-1,i'} \bigl\{}
               - \EE(\bxi^{(i)}) \otimes \EE(\bxi^{(i')}) \otimes \EE(\bxi^{(i)})
               - \EE(\bxi^{(i')}) \otimes [\EE(\bxi^{(i)})]^{\otimes2}\bigr\} \\
  &\quad
     + \sum_{i=1}^p
        Y_{n-1,i}
        \bigl\{\EE[(\bxi^{(i)})^{\otimes3}]
               - \EE[(\bxi^{(i)})^{\otimes2}] \otimes \EE(\bxi^{(i)})
               - \EE(\bxi_{1,1}^{(i)} \otimes \bxi_{1,2}^{(i)} \otimes \bxi_{1,1}^{(i)}) \\
  &\phantom{\quad + \sum_{i=1}^p Y_{n-1,i} \bigl\{}
               - \EE(\bxi^{(i)}) \otimes \EE[(\bxi^{(i)})^{\otimes2}]
               + 2 [\EE(\bxi^{(i)})]^{\otimes3}\bigr\} .
 \end{align*}
Hence
 \begin{equation}\label{Sn1}
  S_{n,1} = \bM_\bxi^{\otimes3} \bY_{n-1}^{\otimes3} + \bA_{3,2}^{(1)} \bY_{n-1}^{\otimes2}
            + \bA_{3,1}^{(1)} \bY_{n-1}
 \end{equation}
 with
 \begin{align*}
  \bA_{3,2}^{(1)}
  &:= \sum_{i=1}^p \sum_{i'=1}^p
       \bigl\{\EE[(\bxi^{(i)})^{\otimes2}] \otimes \EE(\bxi^{(i')})
              + \EE(\bxi^{(i)} \otimes \bxi^{(i')} \otimes \bxi^{(i)})
              + \EE(\bxi^{(i')}) \otimes \EE[(\bxi^{(i)})^{\otimes2}] \\
  &\phantom{:= \sum_{i=1}^p \sum_{i'=1}^p \bigl\{}
              - [\EE(\bxi^{(i)})]^{\otimes2} \otimes \EE(\bxi^{(i')})
              - \EE(\bxi^{(i)}) \otimes \EE(\bxi^{(i')}) \otimes \EE(\bxi^{(i)})
              - \EE(\bxi^{(i')}) \otimes [\EE(\bxi^{(i)})]^{\otimes2}\bigr\} \\
  &\phantom{:= \sum_{i=1}^p \sum_{i'=1}^p}
       \times (\be_i^\top \otimes \be_{i'}^\top)
   \in \RR^{p^3\times p^2} , \\
  \bA_{3,1}^{(1)}
  &:= \sum_{i=1}^p
       \bigl\{\EE[(\bxi^{(i)})^{\otimes3}]
              - \EE[(\bxi^{(i)})^{\otimes2}] \otimes \EE(\bxi^{(i)})
              - \EE(\bxi_{1,1}^{(i)} \otimes \bxi_{1,2}^{(i)} \otimes \bxi_{1,1}^{(i)})
              - \EE(\bxi^{(i)}) \otimes \EE[(\bxi^{(i)})^{\otimes2}] \\
  &\phantom{:= \sum_{i=1}^p \bigl\{}
              + 2 [\EE(\bxi^{(i)})]^{\otimes3}\bigr\}
       \be_i^\top
   \in \RR^{p^3\times p} .
 \end{align*}
Moreover,
 \begin{align*}
  S_{n,2}
  &= \sum_{i=1}^p \sum_{\underset{\SC i'\ne i}{i'=1}}^p
      Y_{n-1,i} Y_{n-1,i'}
      \bigl\{\EE(\bxi^{(i)}) \otimes \EE(\bxi^{(i')}) \otimes \bm_\bvare
             + \EE(\bxi^{(i)}) \otimes \bm_\bvare \otimes \EE(\bxi^{(i')}) \\
  &\phantom{= \sum_{i=1}^p \sum_{\underset{\SC i'\ne i}{i'=1}}^p Y_{n-1,i} Y_{n-1,i'} \bigl\{}
             + \bm_\bvare \otimes \EE(\bxi^{(i)}) \otimes \EE(\bxi^{(i')})\bigr\} \\
  &\quad
     + \sum_{i=1}^p
        Y_{n-1,i} (Y_{n-1,i} - 1)
        \bigl\{[\EE(\bxi^{(i)})]^{\otimes2} \otimes \bm_\bvare
               + \EE(\bxi^{(i)}) \otimes \bm_\bvare \otimes \EE(\bxi^{(i)}) \\
  &\phantom{= + \sum_{i=1}^p Y_{n-1,i} (Y_{n-1,i} - 1) \bigl\{}
               + \bm_\bvare \otimes [\EE(\bxi^{(i)})]^{\otimes2}\bigr\} \\
  &\quad
     + \sum_{i=1}^p
        Y_{n-1,i}
        \bigl\{\EE[(\bxi^{(i)})^{\otimes2}] \otimes \bm_\bvare
               + \EE(\bxi^{(i)} \otimes \bvare \otimes \bxi^{(i)})
               + \bm_\bvare \otimes \EE[(\bxi^{(i)})^{\otimes2}]\bigr\} ,
 \end{align*}
 where \ $\EE(\bxi^{(i)} \otimes \bvare \otimes \bxi^{(i)})$ \ is finite, since there exists a permutation matrix
 \ $\bP \in \RR^{p^2\times p^2}$ \ such that \ $\bu \otimes \bv = \bP (\bv \otimes \bu)$ \ for all
 \ $\bu, \bv \in \RR^p$ \ (see, e.g., Henderson and Searle \cite[formula (6)]{HenSea}),  hence
 \begin{align*}
  \EE(\bxi^{(i)} \otimes \bvare \otimes \bxi^{(i)})
  &= \EE([\bP (\bvare \otimes \bxi^{(i)})] \otimes \bxi^{(i)})
   = \EE\bigl([\bP (\bvare \otimes \bxi^{(i)})] \otimes (\bI_p \bxi^{(i)})\bigr) \\
  &= \EE\bigl((\bP \otimes \bI_p) (\bvare \otimes \bxi^{(i)} \otimes \bxi^{(i)})\bigr)
   = (\bP \otimes \bI_p) \bigl(\bm_\bvare \otimes \EE[(\bxi^{(i)})^{\otimes2}]\bigr) .
 \end{align*}
Thus
 \begin{align*}
  S_{n,2}
  &= \sum_{i=1}^p \sum_{i'=1}^p
      Y_{n-1,i} Y_{n-1,i'}
      \bigl\{\EE(\bxi^{(i)}) \otimes \EE(\bxi^{(i')}) \otimes \bm_\bvare
             + \EE(\bxi^{(i)}) \otimes \bm_\bvare \otimes \EE(\bxi^{(i')}) \\
  &\phantom{= \sum_{i=1}^p \sum_{i'=1}^p Y_{n-1,i} Y_{n-1,i'} \bigl\{}
             + \bm_\bvare \otimes \EE(\bxi^{(i)}) \otimes \EE(\bxi^{(i')})\bigr\} \\
  &\quad
     + \sum_{i=1}^p
        Y_{n-1,i}
        \bigl\{\EE[(\bxi^{(i)})^{\otimes2}] \otimes \bm_\bvare
               + \EE(\bxi^{(i)} \otimes \bvare \otimes \bxi^{(i)})
               + \bm_\bvare \otimes \EE[(\bxi^{(i)})^{\otimes2}] \\
  &\phantom{= + \sum_{i=1}^p Y_{n-1,i} \bigl\{}
               - [\EE(\bxi^{(i)})]^{\otimes2} \otimes \bm_\bvare
               - \EE(\bxi^{(i)}) \otimes \bm_\bvare \otimes \EE(\bxi^{(i)})
               - \bm_\bvare \otimes [\EE(\bxi^{(i)})]^{\otimes2}\bigr\} .
 \end{align*}
Hence
 \begin{equation}\label{Sn2}
  S_{n,2} = \bA_{3,2}^{(2)} \bY_{n-1}^{\otimes2} + \bA_{3,1}^{(2)} \bY_{n-1}
 \end{equation}
 with
 \begin{align*}
  \bA_{3,2}^{(2)}
  &:= \sum_{i=1}^p \sum_{i'=1}^p
       \bigl\{\EE(\bxi^{(i)}) \otimes \EE(\bxi^{(i')}) \otimes \bm_\bvare
              + \EE(\bxi^{(i)}) \otimes \bm_\bvare \otimes \EE(\bxi^{(i')}) \\
  &\phantom{:= \sum_{i=1}^p \sum_{i'=1}^p \bigl\{}
              + \bm_\bvare \otimes \EE(\bxi^{(i)}) \otimes \EE(\bxi^{(i')})\bigr\}
       (\be_i^\top \otimes \be_{i'}^\top)
   \in \RR^{p^3\times p^2} , \\
  \bA_{3,1}^{(2)}
  &:= \sum_{i=1}^p
       \bigl\{\EE[(\bxi^{(i)})^{\otimes2}] \otimes \bm_\bvare
              + \EE(\bxi^{(i)} \otimes \bvare \otimes \bxi^{(i)})
              + \bm_\bvare \otimes \EE[(\bxi^{(i)})^{\otimes2}] \\
  &\phantom{:= \sum_{i=1}^p \bigl\{}
              - [\EE(\bxi^{(i)})]^{\otimes2} \otimes \bm_\bvare
              - \EE(\bxi^{(i)}) \otimes \bm_\bvare \otimes \EE(\bxi^{(i)})
              - \bm_\bvare \otimes [\EE(\bxi^{(i)})]^{\otimes2}\bigr\}
       \be_i^\top
   \in \RR^{p^3\times p} .
 \end{align*}
Further,
 \[
   S_{n,3}
   = \sum_{i=1}^p Y_{n-1,i}
      \bigl\{\EE(\bxi^{(i)}) \otimes \EE(\bvare^{\otimes2})
             + \EE(\bvare \otimes \bxi^{(i)} \otimes \bvare)
             + \EE(\bvare^{\otimes2}) \otimes \EE(\bxi^{(i)})\bigr\}
   = \bA_{3,1}^{(3)} \bY_{n-1}
 \]
 with
 \[
   \bA_{3,1}^{(3)}
   := \sum_{i=1}^p
       \bigl\{\EE(\bxi^{(i)}) \otimes \EE(\bvare^{\otimes2})
              + \EE(\bvare \otimes \bxi^{(i)} \otimes \bvare)
              + \EE(\bvare^{\otimes2}) \otimes \EE(\bxi^{(i)})\bigr\}
       \be_i^\top
   \in \RR^{p^3\times p} ,
 \]
 where \ $\EE(\bvare \otimes \bxi^{(i)} \otimes \bvare)$ \ is finite, since
 \begin{align*}
  \EE(\bvare \otimes \bxi^{(i)} \otimes \bvare)
  &= \EE([\bP (\bxi^{(i)} \otimes \bvare)] \otimes \bvare)
   = \EE\bigl([\bP (\bxi^{(i)} \otimes \bvare)] \otimes (\bI_p \bvare)\bigr) \\
  &= \EE\bigl((\bP \otimes \bI_p) (\bxi^{(i)} \otimes \bvare \otimes \bvare)\bigr)
   = (\bP \otimes \bI_p) \bigl(\EE(\bxi^{(i)}) \otimes \EE[\bvare^{\otimes2}]\bigr) .
 \end{align*}
Consequently, we have
 \[
   \EE(\bY_n^{\otimes3} \mid \cF_{n-1}^\bY)
   = \bM_\bxi^{\otimes3} \bY_{n-1}^{\otimes3} + \bA_{3,2} \bY_{n-1}^{\otimes2} + \bA_{3,1} \bY_{n-1}
     + \EE(\bvare^{\otimes3})
 \]
 with \ $\bA_{3,2} := \bA_{3,2}^{(1)} +\bA_{3,2}^{(2)}$ \ and
 \ $\bA_{3,1} := \bA_{3,1}^{(1)} + \bA_{3,1}^{(2)} + \bA_{3,1}^{(3)}$.
\ Taking the expectation, we get
 \begin{equation}\label{moment_3_Y}
  \EE(\bY_n^{\otimes3})
  = \bM_\bxi^{\otimes3} \EE(\bY_{n-1}^{\otimes3}) + \bA_{3,2} \EE(\bY_{n-1}^{\otimes2}) + \bA_{3,1} \EE(\bY_{n-1})
    + \EE(\bvare^{\otimes3}) .
 \end{equation}
Summarizing, we obtain
 \[
   \begin{bmatrix}
    \EE(\bY_n) \\
    \EE(\bY_n^{\otimes2}) \\
    \EE(\bY_n^{\otimes3})
   \end{bmatrix}
   = \bA_3
     \begin{bmatrix}
      \EE(\bY_{n-1}) \\
      \EE(\bY_{n-1}^{\otimes2}) \\
      \EE(\bY_{n-1}^{\otimes3})
     \end{bmatrix}
     + \begin{bmatrix}
        \bm_\bvare \\
        \EE(\bvare^{\otimes2}) \\
        \EE(\bvare^{\otimes3})
       \end{bmatrix} , \qquad n \in \NN ,
 \]
 with
 \[
   \bA_3 := \begin{bmatrix}
             \bM_\bxi & \bzero & \bzero \\
             \bA_{2,1} & \bM_\bxi^{\otimes2} & \bzero \\
             \bA_{3,1} & \bA_{3,2} & \bM_\bxi^{\otimes3}
            \end{bmatrix}
   \in \RR^{(p+p^2+p^3)\times(p+p^2+p^3)} .
 \]
Taking into account \ $\bY_0 = \bzero$, \ we obtain
 \[
   \begin{bmatrix}
    \EE(\bY_n) \\
    \EE(\bY_n^{\otimes2}) \\
    \EE(\bY_n^{\otimes3})
   \end{bmatrix}
   = \sum_{k=1}^n
      \bA_3^{n-k}
      \begin{bmatrix}
       \bm_\bvare \\
       \EE(\bvare^{\otimes2}) \\
       \EE(\bvare^{\otimes3})
      \end{bmatrix}
   = \sum_{\ell=0}^{n-1}
      \bA_3^\ell
      \begin{bmatrix}
       \bm_\bvare \\
       \EE(\bvare^{\otimes2}) \\
       \EE(\bvare^{\otimes3})
      \end{bmatrix} , \qquad n \in \NN .
 \]
We have \ $\varrho(\bA_3) = \max\{\varrho(\bM_\bxi), \varrho(\bM_\bxi^{\otimes2}), \varrho(\bM_\bxi^{\otimes3})\}$,
 \ where \ $\varrho(\bM_\bxi^{\otimes2}) = [\varrho(\bM_\bxi)]^2$ \ and
 \ \ $\varrho(\bM_\bxi^{\otimes3}) = [\varrho(\bM_\bxi)]^3$.
\ Taking into account \ $\varrho(\bM_\bxi) < 1$, \ we conclude \ $\varrho(\bA_3) = \varrho(\bM_\bxi) < 1$, \ and,
 by \eqref{reprpi}, the third moment of \ $\pi$ \ is finite,
 and
 \begin{equation}\label{moment_3_pi}
  \begin{bmatrix}
    \int_{\RR^p} \bx \, \pi(\dd\bx) \\
    \int_{\RR^p} \bx^{\otimes2} \, \pi(\dd\bx) \\
    \int_{\RR^p} \bx^{\otimes3} \, \pi(\dd\bx) \\
  \end{bmatrix}
  = (\bI_{p+p^2+p^3} - \bA_3)^{-1}
    \begin{bmatrix}
     \bm_\bvare \\
     \EE(\bvare^{\otimes2}) \\
      \EE(\bvare^{\otimes3})
    \end{bmatrix} .
 \end{equation}
Since
 \[
  (\bI_{p+p^2+p^3} - \bA_3)^{-1}
   = \begin{bmatrix}
       (\bI_p - \bM_\bxi)^{-1} & \bzero & \bzero \\
       \bB_{2,1} & (\bI_{p^2} - \bM_\bxi^{\otimes2})^{-1} & \bzero \\
       \bB_{3,1} & \bB_{3,2} & (\bI_{p^3} - \bM_\bxi^{\otimes3})^{-1} \\
     \end{bmatrix},
 \]
 where
 \begin{align*}
  &\bB_{2,1} = (\bI_{p^2} - \bM_\bxi^{\otimes2})^{-1} \bA_{2,1} (\bI_p - \bM_\bxi)^{-1},\\
  &\bB_{3,1} = (\bI_{p^3} - \bM_\bxi^{\otimes3})^{-1}  (\bA_{3,1}(\bI_p - \bM_\bxi)^{-1} + \bA_{3,2}\bB_{2,1}),\\
  &\bB_{3,2} = (\bI_{p^3} - \bM_\bxi^{\otimes3})^{-1} \bA_{3,2} (\bI_{p^2} - \bM_\bxi^{\otimes2})^{-1},
 \end{align*}
 we have
 \[
  \int_{\RR^p} \bx^{\otimes3} \, \pi(\dd\bx)
   = \bB_{3,1} \bm_\bvare + \bB_{3,2} \EE(\bvare^{\otimes2})
     + (\bI_{p^3} - \bM_\bxi^{\otimes3})^{-1} \EE(\bvare^{\otimes3}). \\[-4mm]
 \]
\proofend

\noindent{\bf Proof of Proposition \ref{simple_aggregation}.}
Similarly as \eqref{help12}, we have
 \[
   \EE(\bX_k \mid \cF_{k-1}^\bX) =\bM_\bxi  \bX_{k-1} + \bm_{\bvare}, \qquad k \in \NN ,
 \]
 where \ $\cF_k^\bX := \sigma(\bX_0, \ldots, \bX_k)$, \ $k \in \ZZ_+$.
\ Consequently,
 \begin{align}\label{help3}
  \EE(\bX_k) = \bM_\bxi  \EE(\bX_{k-1}) +  \bm_{\bvare} , \qquad k \in \NN ,
 \end{align}
 and, by \eqref{moment_1_pi},
 \begin{align}\label{help_exp_stat_distr}
  \EE(\bX_0) = (\bI_p - \bM_\bxi )^{-1} \bm_{\bvare} .
 \end{align}
Put
 \begin{align*}
  \bU_k :&= \bX_k - \EE(\bX_k \mid \cF_{k-1}^\bX) =\bX_k - (\bM_\bxi  \bX_{k-1} + \bm_{\bvare}) \\
         &= \sum_{i=1}^p \sum_{\ell=1}^{X_{k-1,i}} (\bxi^{(i)}_{k,\ell} - \EE(\bxi^{(i)}_{k,\ell}))
            + (\bvare_k - \EE(\bvare_k)) ,
          \qquad k \in \NN.
 \end{align*}
Then \ $\EE(\bU_k \mid \cF_{k-1}^\bX) = \bzero$, \ $k \in \NN$, \ and using the independence of
 \ $\bigl\{\bxi^{(i)}_{k,\ell}, \bvare_k : k, \ell \in \NN, i \in \{1,\ldots,p\}\bigr\}$, \ we have
 \begin{equation}\label{help4}
   \EE(U_{k,i}U_{k,j} \mid \cF_{k-1}^\bX)
   = \sum_{q=1}^p  X_{k-1,q} \cov(\xi_{k,1}^{(q,i)}, \xi_{k,1}^{(q,j)}) + \cov(\vare_k^{(i)}, \vare_k^{(j)})
   = \bv_{(i,j)}^\top
      \begin{bmatrix}
       \bX_{k-1} \\
       1
      \end{bmatrix}
 \end{equation}
 for \ $i, j \in \{1, \dots, p\}$ \ and \ $k \in \NN$, \ where \ $[U_{k,1}, \ldots, U_{k,p}]^\top := \bU_k$,
 \ $k \in \NN$.
\ For each \ $k \in \NN$, \ using \ $\bX_k = \bM_\bxi \bX_{k-1} + \bm_{\vare} + \bU_k$ \ and \eqref{help3}, we
 obtain
 \begin{equation}\label{help5}
  \bX_k - \EE(\bX_k) = \bM_\bxi  (\bX_{k-1} - \EE(\bX_{k-1})) + \bU_k , \qquad k \in \NN .
 \end{equation}
Consequently,
 \begin{align*}
  &\EE((\bX_k - \EE(\bX_k)) (\bX_k - \EE(\bX_k))^\top \mid \cF_{k-1}^\bX) \\
  &= \EE((\bM_\bxi (\bX_{k-1} - \EE(\bX_{k-1})) + \bU_k) (\bM_\bxi (\bX_{k-1} - \EE(\bX_{k-1})) + \bU_k)^\top
         \mid \cF_{k-1}^\bX) \\
  &= \EE(\bU_{k} \bU_{k}^\top \mid \cF_{k-1}^\bX)
     + \bM_\bxi (\bX_{k-1} - \EE(\bX_{k-1})) (\bX_{k-1} - \EE(\bX_{k-1}))^\top \bM_\bxi ^\top
 \end{align*}
for all \ $k \in \NN$.
\ Taking the expectation, by \eqref{help_exp_stat_distr} and \eqref{help4}, we conclude
 \[
   \var(\bX_k)
   = \EE(\bU_{k} \bU_{k}^\top) + \bM_\bxi \var(\bX_{k-1}) \bM_\bxi ^\top
   = \bV + \bM_\bxi \var(\bX_{k-1}) \bM_\bxi ^\top ,
   \qquad k \in \NN .
 \]
Under the conditions of the proposition, by Lemma \ref{stationary_moments}, the unique stationary distribution
 \ $\pi$ \ has a finite second moment, hence, using again the stationarity of \ $(\bX_k)_{k\in\ZZ_+}$, \ for each
 \ $N \in \NN$, \ we get
 \begin{align}\label{help10}
  \var(\bX_0)
  = \bV + \bM_\bxi  \var(\bX_0)\bM_\bxi^\top
  = \sum_{k=0}^{N-1} \bM_\bxi^k \bV (\bM_\bxi^\top)^k + \bM_\bxi^N \var(\bX_0) (\bM_\bxi^\top)^N .
 \end{align}
Here \ $\lim_{N\to\infty} \bM_\bxi^N \var(\bX_0) (\bM_\bxi^\top)^N = \bzero \in\RR^{p\times p}$.
\ Indeed, by the Gelfand formula \ $\varrho(\bM_\bxi) = \lim_{k\to\infty} \|\bM_\bxi^k\|^{1/k}$, \ see, e.g.,
 Horn and Johnson \cite[Corollary 5.6.14]{HornJohnson}.
Hence there exists \ $k_0 \in \NN$ \ such that
 \begin{equation}\label{Gelfand}
  \|\bM_\bxi^k\|^{1/k} \leq \varrho(\bM_\bxi) + \frac{1-\varrho(\bM_\bxi)}{2}
  = \frac{1+\varrho(\bM_\bxi)}{2}
  < 1
  \qquad \text{for all \ $k \geq k_0$,}
 \end{equation}
 since \ $\varrho(\bM_\bxi) < 1$.
\ Thus, for all \ $N \geq k_0$,
 \begin{align*}
  \|\bM_\bxi^N \var(\bX_0)(\bM_\bxi^\top)^N\|
  &\leq \|\bM_\bxi^N\| \|\var(\bX_0)\| \|(\bM_\bxi^\top)^N\|
   = \|\bM_\bxi^N\| \|\var(\bX_0)\| \|\bM_\bxi^N\| \\
  &\leq \left(\frac{1+\varrho(\bM_\bxi)}{2}\right)^{2N} \|\var(\bX_0)\| ,
 \end{align*}
 hence \ $\|\bM_\bxi^N \var(\bX_0) (\bM_\bxi^\top)^N\| \to 0$ \ as \ $N \to \infty$.
\ Consequently, \ $\var(\bX_0) =\sum_{k=0}^\infty \bM_\bxi ^{k} \bV (\bM_\bxi ^{\top})^{k}$, \ yielding
 \eqref{var0}.
Moreover, by \eqref{help5},
 \begin{align*}
  &\EE((\bX_0 - \EE(\bX_0)) (\bX_k - \EE(\bX_k))^\top \mid \cF_{k-1}^\bX)
   = (\bX_0 - \EE(\bX_0)) \EE((\bX_k - \EE(\bX_k))^\top \mid \cF_{k-1}^\bX) \\
  &= (\bX_0 - \EE(\bX_0)) (\bX_{k-1} - \EE(\bX_{k-1}) )^\top \bM_\bxi ^\top,
  \qquad k \in \NN .
 \end{align*}
Taking the expectation, we conclude
 \[
   \cov(\bX_0, \bX_k) =  \cov(\bX_0, \bX_{k-1}) \bM_\bxi ^\top , \qquad k \in \NN .
 \]
Hence, by induction, we obtain the formula for \ $\cov(\bX_0, \bX_k)$.
\ The statement will follow from the multidimensional central limit theorem.
Due to the continuous mapping theorem, it is sufficient to show the convergence
 \ $N^{-1/2} (\bS^{(N)}_0, \bS^{(N)}_1, \ldots, \bS^{(N)}_k) \distr (\bcX_0, \bcX_1, \ldots, \bcX_k)$ \ as
 \ $N \to \infty$ \ for all \ $k \in \ZZ_+$.
\ For all \ $k \in \ZZ_+$, \ the random vectors
 \ $\bigl((\bX^{(j)}_0 - \EE(\bX^{(j)}_0))^\top, (\bX^{(j)}_1 - \EE(\bX^{(j)}_1))^\top,
          \ldots, (\bX^{(j)}_k - \EE(\bX^{(j)}_k))^\top\bigr)^\top$,
 \ $j \in \NN$,
 are independent, identically distributed having zero mean vector and covariances
 \[
   \cov(\bX^{(j)}_{\ell_1}, \bX^{(j)}_{\ell_2})
   = \cov(\bX^{(j)}_0, \bX^{(j)}_{\ell_2-\ell_1})
   = \var(\bX_0) (\bM_\bxi ^\top)^{\ell_2-\ell_1}
 \]
 for \ $j \in \NN$, \ $\ell_1, \ell_2 \in \{0, 1, \ldots, k\}$, \ $\ell_1\leq \ell_2$, \ following from the
 strict stationarity of \ $\bX^{(j)}$ \ and from \eqref{cov}.
\proofend


\noindent{\bf Proof of Proposition \ref{simple_aggregation2}.}
It is known that
 \[
   \bU_k
   = \bX_k - \EE(\bX_k \mid \cF^\bX_{k-1})
   = \bX_k - \bM_\bxi  \bX_{k-1} - \bm_{\vare} ,
   \qquad k \in \NN ,
 \]
 are martingale differences with respect to the filtration \ $(\cF_k^\bX)_{k\in\ZZ_+}$.
\ The functional martingale central limit theorem can be applied, see, e.g., Jacod and Shiryaev
 \cite[Theorem VIII.3.33]{JacShi}.
Indeed, using \eqref{help4} and the fact that the first moment of \ $\bX_0$ \ exists and is finite, by
 \eqref{help_ergodic}, for each \ $t \in \RR_+$, \, and \ $i, j \in \{1, \ldots, p\}$, \ we have
 \[
   \frac{1}{n} \sum_{k=1}^\nt \EE(U_{k,i}U_{k,j} \mid \cF^\bX_{k-1})
   \as \bv_{(i,j)}^\top
       \begin{bmatrix}
        \EE(\bX_0) \\
        1
        \end{bmatrix} t
   = V_{i,j} t
   \qquad \text{as \ $n \to \infty$,}
 \]
 and hence the convergence holds in probability as well.
Moreover, the conditional Lindeberg condition holds, namely, for all \ $\delta > 0$,
 \begin{equation}\label{help9}
  \begin{aligned}
   \frac{1}{n}
   \sum_{k=1}^\nt
    \EE\bigl(\|\bU_k\|^2 \bone_{\{\|\bU_k\|>\delta\sqrt{n}\}} \mid \cF^\bX_{k-1}\bigr)
   &\leq \frac{1}{\delta n^{3/2}} \sum_{k=1}^{\nt} \EE(\|\bU_k\|^3 \mid \cF^\bX_{k-1}) \\
   &\leq \frac{C_3 (p+1)^3}{\delta n^{3/2}}
         \sum_{k=1}^\nt \left\|\begin{bmatrix}
                                \bX_{k-1} \\
                                1 \\
                               \end{bmatrix}\right\|^3
    \as 0
  \end{aligned}
 \end{equation}
 with \ $C_3 := \max\{\EE(\|\bxi^{(i)} - \EE(\bxi^{(i)})\|^3)$, \ $i \in \{1, \ldots, p\}$,
 \ $\EE(\|\bvare - \EE(\bvare)\|^3)\}$, \ where the last inequality follows by Proposition 3.3 of Ned\'enyi
 \cite{Ned}, and the almost sure convergence is a consequence of \eqref{help_ergodic}, since, under the
 third order moment assumptions in Proposition \ref{simple_aggregation2}, by Lemma \ref{stationary_moments}
 and {\eqref{help_ergodic},
 \[
   \frac{1}{n}
   \sum_{k=1}^{\nt}
    \left\|\begin{bmatrix}
            \bX_{k-1} \\
            1 \\
           \end{bmatrix}\right\|^3
   \as t \EE\left(\left\|\begin{bmatrix}
                          \bX_0 \\
                          1 \\
                         \end{bmatrix}\right\|^3 \right)
   \qquad \text{as \ $n \to \infty$.}
 \]
Hence we obtain
 \[
   \biggl(\frac{1}{\sqrt{n}} \sum_{k=1}^{\lfloor nt \rfloor} \bU_k\biggr)_{t\in\RR_+}
   \distr \bB \qquad
   \text{as \ $n \to \infty$,}
 \]
 where \ $\bB = (\bB_t)_{t\in\RR_+}$ \ is a \ $p$-dimensional zero mean Brownian motion satisfying
 \ $\var(\bB_1) = \bV$.
\ Using \eqref{help5}, we have
 \[
   \bX_k - \EE(\bX_k)
   = \bM_\bxi ^k (\bX_0 - \EE(\bX_0)) + \sum_{j=1}^k \bM_\bxi ^{k-j} \bU_j ,
   \qquad k \in \NN .
 \]
Consequently, for each \ $n \in \NN$ \ and \ $t \in \RR_+$,
 \begin{equation}\label{help7}
  \begin{aligned}
   &\frac{1}{\sqrt{n}} \sum_{k=1}^\nt (\bX_k - \EE(\bX_k)) \\
   &= \frac{1}{\sqrt{n}}
      \Biggl[\Biggl(\sum_{k=1}^\nt \bM_\bxi^k\Biggr) (\bX_0 - \EE(\bX_0))
             + \sum_{k=1}^\nt \sum_{j=1}^k \bM_\bxi^{k-j} \bU_j\Biggr] \\
  &= \frac{1}{\sqrt{n}}
     \Biggl[(\bI_p - \bM_\bxi)^{-1} (\bM_\bxi - \bM_\bxi^{\lfloor nt \rfloor+1}) (\bX_0 - \EE(\bX_0))
            + \sum_{j=1}^\nt \Biggl(\sum_{k=j}^\nt \bM_\bxi^{k-j}\Biggr) \bU_j \Biggr] \\
  &= \frac{1}{\sqrt{n}}
     \Biggl[(\bI_p - \bM_\bxi)^{-1} (\bM_\bxi - \bM_\bxi^{\lfloor nt \rfloor+1}) (\bX_0 - \EE(\bX_0))
            + (\bI_p - \bM_\bxi)^{-1} \sum_{j=1}^\nt (\bI_p - \bM_\bxi^{\lfloor nt \rfloor-j+1}) \bU_j\Biggr] ,
  \end{aligned}
 \end{equation}
 implying the statement using Slutsky's lemma since \ $\rho(\bM_\bxi) < 1$.
\ Indeed, \ $\lim_{n\to\infty} \bM_\bxi^{\lfloor nt \rfloor+1} = \bzero$ \ by \eqref{Gelfand}, hence
 \[
   \frac{1}{\sqrt{n}}
   (\bI_p - \bM_\bxi)^{-1} (\bM_\bxi - \bM_\bxi^{\lfloor nt \rfloor+1}) (\bX_0 - \EE(\bX_0))
   \as \bzero \qquad \text{as \ $n\to\infty$.}
 \]
Moreover, \ $n^{-1/2} (\bI_p - \bM_\bxi )^{-1} \sum_{j=1}^\nt \bM_\bxi^{\lfloor nt \rfloor-j+1} \bU_j$
 \ converges in \ $L_1$ \ and hence in probability to \ $\bzero$ \ as \ $n \to \infty$, \ since by \eqref{help4},
 \begin{align}\label{help8}
  \EE(|U_{k,j}|)
  \leq \sqrt{\EE(U_{k,j}^2)}
  = \sqrt{\bv_{(j,j)}^\top \begin{bmatrix}\EE(\bX_0)\\ 1\end{bmatrix}}
  = \sqrt{V_{j,j}} ,
  \qquad j \in \{1, \ldots, p\} , \qquad k \in \NN ,
 \end{align}
 and hence
 \begin{align}
  &\EE\biggl(\biggl\|\frac{1}{\sqrt{n}} \sum_{k=1}^\nt \bM_\bxi^{\lfloor nt \rfloor-k+1} \bU_k \biggr\|\biggr)
   \leq \frac{1}{\sqrt{n}} \sum_{k=1}^\nt \EE(\|\bM_\bxi^{\lfloor nt \rfloor-k+1} \bU_k\|) \nonumber \\
  &\leq \frac{1}{\sqrt{n}} \sum_{k=1}^\nt \|\bM_\bxi^{\lfloor nt \rfloor-k+1}\| \EE(\|\bU_k\|)
   \leq \frac{1}{\sqrt{n}} \sum_{k=1}^\nt \|\bM_\bxi^{\lfloor nt \rfloor-k+1}\| \sum_{j=1}^p \EE(|U_{k,j}|)
   \nonumber \\
  &\leq \frac{1}{\sqrt{n}}
        \sum_{k=1}^{\lfloor nt \rfloor} \|\bM_\bxi^{\lfloor nt \rfloor-k+1}\|
         \sum_{j=1}^p \sqrt{V_{j,j}}
   \to 0 \qquad \text{as \ $n \to \infty$,} \label{Gelfand2}
 \end{align}
 since, applying \eqref{Gelfand} for \ $\nt \geq k_0$, \ we have
  \begin{align*}
   &\sum_{k=1}^\nt \|\bM_\bxi^{\lfloor nt \rfloor-k+1}\|
    = \sum_{k=1}^\nt \|\bM_\bxi^k\|
    = \sum_{k=1}^{k_0-1} \|\bM_\bxi^k\| + \sum_{k=k_0}^\nt \|\bM_\bxi^k\| \\
   &\leq \sum_{k=1}^{k_0-1} \|\bM_\bxi^k\|
         + \sum_{k=k_0}^{\lfloor nt \rfloor} \biggl(\frac{1+\varrho(\bM_\bxi)}{2}\biggr)^k
    \leq \sum_{k=1}^{k_0-1} \|\bM_\bxi^k\|
         + \sum_{k=k_0}^\infty \biggl(\frac{1+\varrho(\bM_\bxi)}{2}\biggr)^k
    < \infty .
  \end{align*}
Consequently, by Slutsky's lemma,
 \[
   \biggl(n^{-\frac{1}{2}} \sum_{k=1}^{\lfloor nt \rfloor} (\bX_k - \EE(\bX_k))\biggr)_{t\in\RR_+}
   \distr (\bI_p-\bM_\bxi )^{-1} \bB\,
   \qquad \text{as \ $n \to \infty$,}
 \]
 where \ $\bB = (\bB_t)_{t\in\RR_+}$ \ is a \ $p$-dimensional zero mean Brownian motion satisfying
 \ $\var(\bB_1) = \bV$, \ as desired.
\proofend


\noindent{\bf Proof of Theorem \ref{double_aggregation}.}
First, we prove \eqref{help2_double_aggregation}.
For all \ $N, m \in \NN$ \ and all \ $t_1, \ldots, t_m \in\RR_+$, \ by Proposition \ref{simple_aggregation2} and
 the continuity theorem, we have
 \[
   \frac{1}{\sqrt{n}} (\bS^{(N,n)}_{t_1}, \ldots, \bS^{(N,n)}_{t_m})
   \distr (\bI_p-\bM_\bxi)^{-1} \sum_{\ell=1}^N (\bB^{(\ell)}_{t_1}, \ldots, \bB^{(\ell)}_{t_m})
 \]
 as \ $n \to \infty$, \ where \ $\bB^{(\ell)} = (\bB^{(\ell)}_t)_{t\in\RR_+}$, \ $\ell \in \{1, \ldots, N\}$, \ are
 independent \ $p$-dimensional zero mean Brownian motions satisfying \ $\var(\bB^{(\ell)}_1) = \bV$,
 \ $\ell \in \{1, \ldots, p\}$.
\ Since
 \[
   \frac{1}{\sqrt{N}} \sum_{\ell=1}^N (\bB^{(\ell)}_{t_1}, \ldots, \bB^{(\ell)}_{t_m})
   \distre (\bB_{t_1}, \ldots, \bB_{t_m}) , \qquad N \in \NN , \quad  m \in \NN ,
 \]
 we obtain the convergence \eqref{help2_double_aggregation}.

Now, we turn to prove \eqref{help1_double_aggregation}.
For all \ $n \in \NN$ \ and for all \ $t_1, \ldots, t_m \in \RR_+$ \ with \ $t_1 < \ldots < t_m$, \ $m \in \NN$,
 \ by Proposition \ref{simple_aggregation} and by the continuous mapping theorem, we have
 \begin{align*}
  \frac{1}{\sqrt{N}} \bigl((\bS^{(N,n)}_{t_1})^\top , \ldots, (\bS^{(N,n)}_{t_m})^\top\bigr)^ \top
  &\distr \Biggl(\sum_{k=1}^{\lfloor nt_1\rfloor} \bcX_k^\top , \ldots,
                 \sum_{k=1}^{\lfloor nt_m\rfloor} \bcX_k^\top\Biggr)^\top \\
  &\distre \cN_ {pm} \Biggl(\bzero,
                            \var\Biggl(\Biggl(\sum_{k=1}^{\lfloor nt_1\rfloor}
                                               \bcX_k^\top , \ldots,
                                               \sum_{k=1}^{\lfloor nt_m\rfloor}
                                                \bcX_k^\top\Biggr)^\top\Biggr)\Biggr)
 \end{align*}
 as \ $N \to \infty$, \ where \ $(\bcX_k)_{k\in\ZZ_+}$ \ is the \ $p$-dimensional zero mean stationary
 Gaussian process given in Proposition \ref{simple_aggregation} and, by \eqref{cov},
 \begin{align*}
  &\var\Biggl(\Biggl(\sum_{k=1}^{\lfloor nt_1\rfloor} \bcX_k^\top , \ldots,
                     \sum_{k=1}^{\lfloor nt_m\rfloor} \bcX_k^\top\Biggr)^\top\Biggr)
   = \left(\cov\Biggl(\sum_{k=1}^{\lfloor nt_i\rfloor} \bcX_k,
                      \sum_{k=1}^{\lfloor nt_j\rfloor} \bcX_k\Biggr) \right)_{i,j=1}^m
 \end{align*}
 \begin{align*}
  &= \left(\sum_{k=1}^{\lfloor nt_i\rfloor} \sum_{\ell=1}^{\lfloor nt_j\rfloor}
            \cov(\bcX_k, \bcX_\ell)\right)_{i,j=1}^m \\
  &= \Bigg(\sum_{k=1}^{\lfloor nt_i\rfloor} \sum_{\ell=1}^{(k-1)\wedge \lfloor nt_j\rfloor}
            \bM_\bxi ^{k-\ell} \var(\bX_0)
           + {(\lfloor nt_i\rfloor\wedge \lfloor nt_j\rfloor)}\var(\bX_0)\\
&\phantom{= \Bigg(\;}
           + \var(\bX_0)
             \sum_{k=1}^{\lfloor nt_i\rfloor} \sum_{\ell=k+1}^{\lfloor nt_j\rfloor}
              (\bM_\bxi^\top)^{\ell-k}\Bigg)_{i,j=1}^m ,
 \end{align*}
 where \ $\sum_{\ell=q_1}^{q_2} := 0$ \ for all \ $q_2 < q_1$, \ $q_1, q_2 \in \NN$.
\ By the continuity theorem, for all \ $\btheta_1, \ldots, \btheta_m \in \RR^p$, \ $m \in \NN$, \ we conclude
 \begin{align*}
  &\lim_{N\to\infty}
    \EE\biggl(\exp\biggl\{\ii
                          \sum_{j=1}^m
                           \btheta_j^\top n^{-1/2} N^{-1/2}
                           \bS^{(N,n)}_{t_j}\biggr\}\biggr) \\
  &= \exp\left\{-\frac{1}{2n}
                  \sum_{i=1}^m \sum_{j=1}^m
                   \btheta_i^\top \left[
                   \sum_{k=1}^{\lfloor nt_i\rfloor}
                    \sum_{\ell=1}^{\lfloor nt_j\rfloor}
                    \cov(\bcX_k,\bcX_\ell)\right]
			\btheta_j \right\}\\
  &\to  \exp\biggl\{-\frac{1}{2}
                    \sum_{i=1}^m \sum_{j=1}^m
                     (t_i \land t_j) \btheta_i^\top \Big[
 \bM_\bxi (\bI_p-\bM_\bxi )^{-1}\var(\bX_0)+\var(\bX_0) \\
&\phantom{\to  \exp\biggl\{-\frac{1}{2} \sum_{i=1}^m \sum_{j=1}^m (t_i \land t_j) \btheta_i^\top  \Big[\,}
  +\var(\bX_0)
	(\bI_p-\bM_\bxi ^\top)^{-1}\bM_\bxi ^\top\Big] \btheta_j  \biggr\}
  \qquad \text{as \ $n \to \infty$.}
 \end{align*}
Indeed, for all \ $s, t \in \RR_+$ \ with \ $s < t$, \ we have
  \begin{align*}
   &\frac{1}{n} \sum_{k=1}^\ns \sum_{\ell=1}^\nt \cov(\bcX_k, \bcX_\ell) \\
   &= \frac{1}{n} \sum_{k=1}^\ns \sum_{\ell=1}^{k-1} \bM_\bxi^{k-\ell} \var(\bX_0)
      + \frac{\ns}{n}\var(\bX_0)
      + \frac{1}{n} \var(\bX_0)\sum_{k=1}^\ns \sum_{\ell=k+1}^\nt (\bM_\bxi^\top)^{\ell-k} \\
   &= \frac{1}{n} \sum_{k=1}^\ns (\bM_\bxi - \bM_\bxi^k) (\bI_p - \bM_\bxi)^{-1} \var(\bX_0)
      + \frac{\ns}{n}\var(\bX_0) \\
   &\phantom{=\,}
      + \frac{1}{n} \var(\bX_0) (\bI_p - \bM_\bxi^\top)^{-1}
        \sum_{k=1}^\ns (\bM_\bxi^\top - (\bM_\bxi^\top)^{\lfloor nt\rfloor-k+1}) \\
   &= \frac{1}{n}
      \Bigl(\ns \bM_\bxi - \bM_\bxi (\bI_p - \bM_\bxi^\ns) (\bI_p-\bM_\bxi )^{-1}\Bigr)
      (\bI_p - \bM_\bxi)^{-1} \var(\bX_0)
      + \frac{\lfloor ns\rfloor}{n} \var(\bX_0) \\
   &\quad
      + \frac{1}{n} \var(\bX_0) (\bI_p - \bM_\bxi^\top)^{-1}
        \Bigl(\ns \bM_\bxi^\top
              - (\bI_p - \bM_\bxi^\top)^{-1}
                (\bI_p - (\bM_\bxi^\top)^\ns) (\bM_\bxi ^\top)^{\nt - \ns + 1}\Bigr)
  \end{align*}
  \begin{align*}
   &= \frac{\ns}{n}
      \Bigl(\bM_\bxi (\bI_p - \bM_\bxi)^{-1} \var(\bX_0) + \var(\bX_0)
            + \var(\bX_0) (\bI_p - \bM_\bxi^\top)^{-1} \bM_\bxi ^\top\Bigr) \\
   &\quad
      - \frac{1}{n} \Bigl(\bM_\bxi (\bI_p - \bM_\bxi^\ns) (\bI_p - \bM_\bxi)^{-2} \var(\bX_0) \\
   &\phantom{\quad-\frac{1}{n}\Big[\;}
      + \var(\bX_0) (\bI_p - \bM_\bxi ^\top)^{-2}
        (\bI_p - (\bM_\bxi^\top)^\ns) (\bM_\bxi^\top)^{\nt-\ns+1}\Bigr) \\
   &\to s \Bigl(\bM_\bxi (\bI_p - \bM_\bxi)^{-1} \var(\bX_0) + \var(\bX_0)
                + \var(\bX_0) (\bI_p - \bM_\bxi^\top)^{-1} \bM_\bxi^\top\Bigr) \qquad \text{as \ $n \to \infty$,}
  \end{align*}
 since \ $\lim_{n\to\infty} \bM_\bxi^\ns = \bzero$, \ $\lim_{n\to\infty} (\bM_\bxi^\top)^\ns = \bzero$ \ and
 \ $\lim_{n\to\infty} (\bM_\bxi^\top)^{\nt-\ns+1} = \bzero$ \ by \eqref{Gelfand}.
It remains to show that
 \begin{equation}\label{help6}
  \begin{aligned}
   &\bM_\bxi (\bI_p - \bM_\bxi)^{-1} \var(\bX_0) + \var(\bX_0)
    + \var(\bX_0) (\bI_p - \bM_\bxi^\top)^{-1} \bM_\bxi^\top \\
   &= (\bI_p - \bM_\bxi)^{-1} \bV (\bI_p - \bM_\bxi^\top)^{-1} .
  \end{aligned}
 \end{equation}
We have
 \begin{align}\label{help11}
  \bM_\bxi (\bI_p - \bM_\bxi)^{-1}
  = (\bI_p - (\bI_p - \bM_\bxi)) (\bI_p - \bM_\bxi)^{-1}
  = (\bI_p - \bM_\bxi)^{-1} - \bI_p ,
 \end{align}
 and hence \ $(\bI_p - \bM_\bxi^\top)^{-1} \bM_\bxi^\top = (\bI_p - \bM_\bxi^\top)^{-1} - \bI_p$, \ thus the
 left-hand side of equation \eqref{help6} can be written as
 \begin{align*}
  &((\bI_p - \bM_\bxi)^{-1} - \bI_p) \var(\bX_0) + \var(\bX_0)
   + \var(\bX_0) ((\bI_p - \bM_\bxi^\top)^{-1} - \bI_p) \\
  &= (\bI_p - \bM_\bxi)^{-1} \var(\bX_0) - \var(\bX_0) + \var(\bX_0) (\bI_p - \bM_\bxi^\top)^{-1} .
 \end{align*}
By \eqref{help10}, we have \ $\bV = \var(\bX_0) - \bM_\bxi \var(\bX_0) \bM_\bxi^\top$, \ hence, by \eqref{help11},
 the right-hand side of the equation \eqref{help6} can be written as
 \begin{align*}
  &(\bI_p - \bM_\bxi)^{-1} (\var(\bX_0) - \bM_\bxi \var(\bX_0) \bM_\bxi^\top) (\bI_p - \bM_\bxi^\top)^{-1} \\
  &= (\bI_p - \bM_\bxi)^{-1} \var(\bX_0) (\bI_p - \bM_\bxi^\top)^{-1}
     - (\bI_p - \bM_\bxi)^{-1} \bM_\bxi \var(\bX_0) \bM_\bxi^\top (\bI_p - \bM_\bxi^\top)^{-1} \\
  &= (\bI_p - \bM_\bxi)^{-1} \var(\bX_0) (\bI_p - \bM_\bxi^\top)^{-1}
     - ((\bI_p - \bM_\bxi)^{-1} - \bI_p) \var(\bX_0) ((\bI_p - \bM_\bxi^\top)^{-1} - \bI_p) \\
  &= (\bI_p - \bM_\bxi)^{-1} \var(\bX_0) - \var(\bX_0) + \var(\bX_0) (\bI_p - \bM_\bxi^\top)^{-1} ,
 \end{align*}
 and we conclude \eqref{help6}.
This implies the convergence \eqref{help1_double_aggregation}.
\proofend


\noindent{\bf Proof of Theorem \ref{simulataneous_aggregation}.}
As \ $n$ \ and \ $ N$ \ converge to infinity simultaneously, \eqref{help3_simulataneous_aggregation} is equivalent
 to \ $(nN_n)^{-\frac{1}{2}} \bS^{(N_n,n)} \distr (\bI_p-\bM_\bxi )^{-1} \, \bB$ \ as \ $n \to \infty$ \ for any
 sequence \ $(N_n)_{n\in\NN}$ \ of positive integers such that \ $\lim_{n\to\infty} N_n = \infty$.
\ As we have seen in the proof of Proposition \ref{simple_aggregation2}, for each \ $j \in \NN$,
 \[
   \bU_k^{(j)}
   := \bX_k^{(j)} - \EE(\bX_k^{(j)} \mid \cF^{\bX_{k-1}^{(j)}})
   = \bX_k^{(j)} - \bM_\bxi  \bX_{k-1}^{(j)} - \bm_{\vare} ,
   \qquad k \in \NN ,
 \]
 are martingale differences with respect to the filtration \ $(\cF_k^{\bX^{(j)}})_{k\in\ZZ_+}$.
\ We are going to apply the functional martingale central limit theorem, see, e.g., Jacod and Shiryaev
 \cite[Theorem VIII.3.33]{JacShi}, for the triangular array consisting of the random vectors
 \[
   (\bV_k^{(n)})_{k\in\NN}
   := (nN_n)^{-\frac{1}{2}}
      \bigl(\bU_1^{(1)}, \ldots, \bU_1^{(N_n)}, \bU_2^{(1)}, \ldots, \bU_2^{(N_n)} ,
            \bU_3^{(1)}, \ldots, \bU_3^{(N_n)}, \ldots\bigr)
 \]
 in the \ $n^\mathrm{th}$ \ row for each \ $n \in \NN$ \ with the filtration \ $(\cF_k^{(n)})_{k\in\ZZ_+}$ \ given
 by \ $\cF_k^{(n)} := \cF_k^{\bY^{(n)}} = \sigma(\bY_0^{(n)}, \ldots, \bY_k^{(n)})$, \ where
 \[
   (\bY_k^{(n)})_{k\in\ZZ_+}
   := \bigl((\bX_0^{(1)}, \ldots, \bX_0^{(N_n)}), \bX_1^{(1)}, \ldots, \bX_1^{(N_n)},
             \bX_2^{(1)}, \ldots, \bX_2^{(N_n)}, \ldots\bigr) .
 \]
Hence \ $\cF_0^{(n)} = \sigma(\bX_0^{(1)}, \ldots, \bX_0^{(N_n)})$, \ and for each \ $k = \ell N_n + r$ \ with
 \ $\ell \in \ZZ_+$ \ and \ $r \in \{1, \ldots, N_n\}$, \ we have
 \[
   \cF_k^{(n)} = \sigma\bigl(\bigl(\cup_{j=1}^r \cF_{\ell+1}^{\bX^{(j)}}\bigr)
                             \cup \bigl(\cup_{j=r+1}^{N_n} \cF_\ell^{\bX^{(j)}}\bigr)\bigr) ,
 \]
 where \ $\cup_{j=N_n+1}^{N_n} := \emptyset$.
\ Moreover, \ $\bY_0^{(n)} = (\bX_0^{(1)}, \ldots, \bX_0^{(N_n)})$, \ and for \ $k = \ell N_n + r$ \ with
 \ $\ell \in \ZZ_+$ \ and \ $r \in \{1, \ldots, N_n\}$, \ we have \ $\bY_k^{(n)} = \bX_{\ell+1}^{(r)}$ \ and
 \ $\bV_k^{(n)} = (nN_n)^{-\frac{1}{2}} \bU_{\ell+1}^{(r)}$.

Next we check that for each \ $n \in \NN$, \ $(\bV_k^{(n)})_{k\in\NN}$ \ is a sequence of martingale differences
 with respect to \ $(\cF_k^{(n)})_{k\in\ZZ_+}$.
\ We will use that \ $\EE(\bxi \mid \sigma(\cG_1 \cup \cG_2)) = \EE(\bxi \mid \cG_1)$ \ for a random vector
 \ $\bxi$ \ and for \ $\sigma$-algebras \ $\cG_1 \subset \cF$ \ and \ $\cG_2 \subset \cF$ \ such that
 \ $\sigma(\sigma(\bxi) \cup \cG_1)$ \ and \ $\cG_2$ \ are independent and \ $\EE(\|\bxi\|) < \infty$.
\ For each \ $k = \ell N_n + 1$ \ with \ $\ell \in \ZZ_+$, \ we have
 \ $\EE(\bV_k^{(n)} \mid \cF_{k-1}^{(n)}) = (nN_n)^{-\frac{1}{2}} \EE(\bU_{\ell+1}^{(1)} \mid \cF_\ell^{\bX^{(1)}})
    = \bzero$, \ since
 \[
   \EE(\bU_{\ell+1}^{(1)} \mid \cF_{k-1}^{(n)})
   = \EE(\bU_{\ell+1}^{(1)} \mid \sigma(\cup_{j=1}^{N_n} \cF_\ell^{\bX^{(j)}}))
   = \EE(\bU_{\ell+1}^{(1)} \mid \cF_\ell^{\bX^{(1)}})
   = \bzero .
 \]
In a similar way, for each \ $k = \ell N_n + r$ \ with \ $\ell \in \ZZ_+$ \ and \ $r \in \{2, \ldots, N_n\}$, \ we
 have
 \ $\EE(\bV_k^{(n)} \mid \cF_{k-1}^{(n)})
    = (nN_n)^{-\frac{1}{2}} \EE(\bU_{\ell+1}^{(r)} \mid \cF_\ell^{\bX^{(r)}}) = \bzero$,
 \ since
 \[
   \EE(\bU_{\ell+1}^{(r)} \mid \cF_{k-1}^{(n)})
   = \EE(\bU_{\ell+1}^{(r)}
         \mid \sigma((\cup_{j=1}^{r-1} \cF_{\ell+1}^{\bX^{(j)}}) \cup (\cup_{j=r}^{N_n} \cF_\ell^{\bX^{(j)}}))) \\
   = \EE(\bU_{\ell+1}^{(r)} \mid \cF_\ell^{\bX^{(r)}})
   = \bzero .
 \]
We want to obtain a functional central limit theorem for the sequence
 \[
   \Biggl(\sum_{k=1}^{\lfloor nt \rfloor N_n} \bV_k^{(n)}\Biggr)_{t\in\RR_+}
   = \biggl(\frac{1}{\sqrt{nN_n}} \sum_{\ell=1}^\nt \sum_{r=1}^{N_n} \bU_\ell^{(r)}\biggr)_{t\in\RR_+} , \qquad
   n \in \NN .
 \]
First, we calculate the conditional variance matrix of \ $\bV_k^{(n)}$.
\ If \ $k = \ell N_n + 1$ \ with \ $\ell \in \ZZ_+$, \ then
 \begin{align*}
  \EE(\bV_k^{(n)} (\bV_k^{(n)})^\top \mid \cF_{k-1}^{(n)})
  &= (nN_n)^{-1}
     \EE(\bU_{\ell+1}^{(1)} (\bU_{\ell+1}^{(1)})^\top \mid \sigma(\cup_{j=1}^{N_n} \cF_\ell^{\bX^{(j)}})) \\
  &= (nN_n)^{-1} \EE(\bU_{\ell+1}^{(1)} (\bU_{\ell+1}^{(1)})^\top \mid \cF_\ell^{\bX^{(1)}}) .
 \end{align*}
In a similar way, if \ $k = \ell N_n + r$ \ with \ $\ell \in \ZZ_+$ \ and \ $r \in \{2, \ldots, N_n\}$, \ then
 \begin{align*}
  \EE(\bV_k^{(n)} (\bV_k^{(n)})^\top \mid \cF_{k-1}^{(n)})
  &= (nN_n)^{-1}
     \EE(\bU_{\ell+1}^{(r)} (\bU_{\ell+1}^{(r)})^\top
         \mid \sigma((\cup_{j=1}^{r-1} \cF_{\ell+1}^{\bX^{(j)}}) \cup (\cup_{j=r}^{N_n} \cF_\ell^{\bX^{(j)}}))) \\
  &= (nN_n)^{-1} \EE(\bU_{\ell+1}^{(r)} (\bU_{\ell+1}^{(r)})^\top \mid \cF_\ell^{\bX^{(r)}}) .
 \end{align*}
Consequently, for each \ $n \in \NN$ \ and \ $t \in \RR_+$, \ we have
 \begin{align*}
  \sum_{k=1}^{\nt N_n} \EE(\bV_k^{(n)} (\bV_k^{(n)})^\top \mid \cF_{k-1}^{(n)})
  &= \sum_{\ell=1}^\nt \sum_{r=1}^{N_n}
      \EE(\bV_{(\ell-1)N_n+r}^{(n)} (\bV_{(\ell-1)N_n+r}^{(n)})^\top \mid \cF_{(\ell-1)N_n+r-1}^{(n)}) \\
  &= \frac{1}{nN_n}
     \sum_{\ell=1}^\nt \sum_{r=1}^{N_n}
      \EE(\bU_\ell^{(r)} (\bU_\ell^{(r)})^\top \mid \cF_{\ell-1}^{\bX^{(r)}}) .
 \end{align*}
Next, we show that for each \ $t \in \RR_+$ \ and \ $i, j \in \{1,\ldots,p\}$, \ we have
 \[
   \frac{1}{nN_n}
   \sum_{\ell=1}^\nt \sum_{r=1}^{N_n} \EE(U^{(r)}_{\ell,i}U^{(r)}_{\ell,j} \mid \cF^{\bX^{(r)}}_{\ell-1})
   = \frac{1}{nN_n}
     \sum_{\ell=1}^\nt \sum_{r=1}^{N_n}
      \bv_{(i,j)}^\top
      \begin{bmatrix}
       \bX_{\ell-1}^{(r)} \\
       1
      \end{bmatrix}
   \stoch
   \bv_{(i,j)}^\top
   \begin{bmatrix}
    \EE(\bX_0) \\
    1
   \end{bmatrix}
   t
   = V_{i,j} t
 \]
 as \ $n \to \infty$.
\ Indeed, the equality follows by \eqref{help4}, and for the convergence in probability, note that
 \ $\lim_{n\to\infty} \frac{\nt}{n} = t$, \ $t \in \RR_+$, \ and, by Cauchy-Schwarz inequality,
 \begin{align*}
  &\EE\left(\left(\frac{1}{\nt N_n}
                  \sum_{\ell=1}^{\nt} \sum_{r=1}^{N_n} \bv_{(i,j)}^\top
                   \begin{bmatrix} \bX_{\ell-1}^{(r)} - \EE(\bX_0) \\ 0 \end{bmatrix}\right)^2\right) \\
  &= \frac{1}{\nt^2 N_n^2}
     \EE\left(\left(\bv_{(i,j)}^\top
                    \sum_{\ell_1=1}^\nt \sum_{r_1=1}^{N_n}
                     \begin{bmatrix} \bX_{\ell_1-1}^{(r_1)} - \EE(\bX_0) \\ 0 \end{bmatrix}\right)
              \left(\sum_{\ell_2=1}^\nt \sum_{r_2=1}^{N_n}
                    \begin{bmatrix} \bX_{\ell_2-1}^{(r_2)} - \EE(\bX_0) \\ 0 \end{bmatrix}^\top
              \bv_{(i,j)}\right)\right) \\
  &= \frac{1}{\nt^2 N_n^2} \bv_{(i,j)}^\top
     \sum_{\ell_1=1}^\nt \sum_{\ell_2=1}^\nt \sum_{r_1=1}^{N_n} \sum_{r_2=1}^{N_n}
      \begin{bmatrix}
       \EE((\bX_{\ell_1-1}^{(r_1)} - \EE(\bX_0)) (\bX_{\ell_2-1}^{(r_2)} - \EE(\bX_0))^\top) & \bzero \\
       \bzero & 0 \\
      \end{bmatrix}
     \bv_{(i,j)} \\
  &= \frac{1}{\nt^2 N_n} \bv_{(i,j)}^\top
     \sum_{\ell_1=1}^\nt \sum_{\ell_2=1}^\nt
      \begin{bmatrix}
       \EE((\bX_{\ell_1-1} - \EE(\bX_0))(\bX_{\ell_2-1} - \EE(\bX_0))^\top) & \bzero \\
       \bzero & 0 \\
      \end{bmatrix}
      \bv_{(i,j)} \\
  &\leq \frac{1}{\nt^2 N_n} \|\bv_{(i,j)}\|^2
        \sum_{\ell_1=1}^\nt \sum_{\ell_2=1}^\nt
         \EE\big(\|(\bX_{\ell_1-1} - \EE(\bX_0)) (\bX_{\ell_2-1} - \EE(\bX_0))^\top\|\big) \\
  &\leq \frac{1}{\nt^2 N_n} \|\bv_{(i,j)}\|^2
        \sum_{\ell_1=1}^\nt \sum_{\ell_2=1}^\nt \sum_{m_1=1}^p  \sum_{m_2=1}^p
         \EE(|(X_{\ell_1-1,m_1} - \EE(X_{0,m_1})) (X_{\ell_2-1,m_2} - \EE(X_{0,m_2}))|) \\
  &\leq \frac{1}{\nt^2 N_n} \|\bv_{(i,j)}\|^2
        \sum_{\ell_1=1}^\nt \sum_{\ell_2=1}^\nt
        \sum_{m_1=1}^{p}\sum_{m_2=1}^p
         \sqrt{\var(X_{\ell_1-1,m_1}) \var(X_{\ell_2-1,m_2})} \\
  &= \frac{1}{N_n} \|\bv_{(i,j)}\|^2
     \sum_{m_1=1}^p \sum_{m_2=1}^p \sqrt{\var(X_{0,m_1}) \var(X_{0,m_2})}
   \to 0 \qquad \text{as \ $n \to \infty$,}
 \end{align*}
 where we used that \ $\|\bQ\| \leq \sum_{i=1}^p \sum_{j=1}^p |q_{i,j}|$ \ for every matrix
 \ $\bQ = (q_{i,j})_{i,j=1}^p \in \RR^{p\times p}$.

Moreover, in a similar way, the conditional Lindeberg condition holds, namely, for all \ $\delta > 0$,
 \begin{align*}
  \sum_{k=1}^{\nt N_n} \EE(\|\bV_k^{(n)}\|^2 \bbone_{\{\|\bV_k^{(n)}\|>\delta\}} \mid \cF_{k-1}^{(n)})
  &= \frac{1}{nN_n}
     \sum_{\ell=1}^\nt \sum_{r=1}^{N_n}
      \EE(\|\bU_\ell^{(r)}\|^2 \bbone_{\{\|\bU_\ell^{(r)}\|>\delta\sqrt{nN_n}\}} \mid \cF^{\bX^{(r)}}_{\ell-1}) \\
  &\leq \frac{1}{\delta n^{3/2} N_n^{1/2}}
        \sum_{\ell=1}^\nt \EE(\|\bU_\ell^{(1)}\|^3 \mid \cF^{\bX^{(1)}}_{\ell-1})
   \as 0 \qquad\text{as \ $n \to \infty$,}
 \end{align*}
 where the almost sure convergence follows by \eqref{help9}.
Hence we obtain
 \[
   \biggl(\frac{1}{\sqrt{nN_n}} \sum_{\ell=1}^\nt \sum_{r=1}^{N_n} \bU_\ell^{(r)}\biggr)_{t\in\RR_+}
   = \Biggl(\sum_{k=1}^{\lfloor nt \rfloor N_n} \bV_k^{(n)}\Biggr)_{t\in\RR_+}
   \distr \bB \qquad
   \text{as \ $n \to \infty$,}
 \]
 where \ $\bB = (\bB_t)_{t\in\RR_+}$ \ is a \ $p$-dimensional zero mean Brownian motion satisfying
 \ $\var(\bB_1) = \bV$.
\ Using \eqref{help7}, for each \ $n \in \NN$ \ and \ $t \in \RR_+$, \ we have
 \begin{align*}
  &\frac{1}{\sqrt{nN_n}}
   \sum_{\ell=1}^\nt \sum_{r=1}^{N_n} (\bX_\ell^{(r)} - \EE(\bX_\ell^{(r)})) \\
  &= \frac{1}{\sqrt{n}}
     \Bigg[(\bI_p - \bM_\bxi )^{-1} (\bM_\bxi - \bM_\bxi ^{\lfloor nt \rfloor+1})
     \frac{1}{\sqrt{N_n}} \sum_{r=1}^{N_n} (\bX_0^{(r)} - \EE(\bX_0^{(r)}))\Bigg] \\
  &\quad
     - \frac{1}{\sqrt{n}}
       \Bigg[(\bI_p - \bM_\bxi)^{-1}
             \sum_{m=1}^\nt
              \bM_\bxi ^{\lfloor nt \rfloor-m+1} \frac{1}{\sqrt{N_n}} \sum_{r=1}^{N_n} \bU_m^{(r)}\Bigg]
       + (\bI_p - \bM_\bxi)^{-1} \frac{1}{\sqrt{nN_n}}
         \sum_{m=1}^\nt \sum_{r=1}^{N_n} \bU_m^{(r)} ,
 \end{align*}
 implying the statement using Slutsky's lemma, since \ $\rho(\bM_\bxi) < 1$.
\ Indeed, \ $\lim_{n\to\infty} \bM_\bxi ^{\lfloor nt \rfloor+1} = \bzero$ \ by \eqref{Gelfand}, thus
 \[
   \lim_{n\to\infty} (\bI_p-\bM_\bxi )^{-1} (\bM_\bxi - \bM_\bxi ^{\lfloor nt \rfloor+1})
   = (\bI_p - \bM_\bxi )^{-1} \bM_\bxi ,
 \]
 and, by Proposition \ref{simple_aggregation},
 \[
   \frac{1}{\sqrt{N_n}} \sum_{r=1}^{N_n} (\bX_0^{(r)} - \EE(\bX_0^{(r)}))
   \distr \cN_p(\bzero, \var(\bX_0))
   \qquad \text{as \ $n\to\infty$,}
 \]
 where \ $\cN_p(\bzero, \var(\bX_0))$ \ denotes a \ $p$-dimensional normal distribution with zero mean and with
 covariance matrix \ $\var(\bX_0)$, \ and then Slutsky's lemma yields that
 \[
   \frac{1}{\sqrt{n}}
   \Bigg[(\bI_p - \bM_\bxi )^{-1} (\bM_\bxi -\bM_\bxi ^{\lfloor nt \rfloor+1})
         \frac{1}{\sqrt{N_n}} \sum_{r=1}^{N_n} (\bX_0^{(r)} - \EE(\bX_0^{(r)}))\Bigg]
   \stoch \bzero \qquad \text{as \ $n \to \infty$.}
 \]
Further,
 \begin{align*}
  &\biggl\|\EE\biggl(\frac{1}{\sqrt{n}}
                     \sum_{m=1}^\nt
                      \bM_\bxi^{\lfloor nt \rfloor-m+1} \frac{1}{\sqrt{N_n}}
                      \sum_{r=1}^{N_n}\bU_m^{(r)}\biggr)\biggr\|
   \leq \frac{1}{\sqrt{n}}
        \sum_{m=1}^\nt
         \EE\biggl(\biggl\|\bM_\bxi^{\lfloor nt \rfloor-m+1} \frac{1}{\sqrt{N_n}}
                           \sum_{r=1}^{N_n} \bU_m^{(r)}\biggr\|\biggr) \\
  &\leq \frac{1}{\sqrt{n}}
        \sum_{m=1}^\nt
         \|\bM_\bxi^{\lfloor nt \rfloor-m+1}\| \EE\biggl(\biggl\|\frac{1}{\sqrt{N_n}}
         \sum_{r=1}^{N_n} \bU_m^{(r)}\biggr\|\bigg) \\
   &\leq \frac{1}{\sqrt{n}}
         \sum_{m=1}^\nt
          \|\bM_\bxi^{\lfloor nt \rfloor-m+1}\|
          \sum_{j=1}^p
           \EE\biggl(\biggl|\frac{1}{\sqrt{N_n}} \sum_{r=1}^{N_n} U^{(r)}_{m,j} \biggr|\biggr)
 \end{align*}
 \begin{align*}
   &\leq \frac{1}{\sqrt{n}}
         \sum_{m=1}^\nt
          \|\bM_\bxi^{\lfloor nt \rfloor-m+1}\|
          \sum_{j=1}^p
           \sqrt{\EE\biggl(\biggl(\frac{1}{\sqrt{N_n}} \sum_{r=1}^{N_n} U^{(r)}_{m,j}\biggr)^2\biggr)} \\
   &= \frac{1}{\sqrt{n}}
      \sum_{m=1}^\nt \|\bM_\bxi^{\lfloor nt \rfloor-m+1}\| \sum_{j=1}^p \sqrt{\EE((U^{(1)}_{m,j})^2)} \\
   &\leq \frac{1}{\sqrt{n}}
         \sum_{m=1}^\nt
          \|\bM_\bxi^{\lfloor nt \rfloor-m+1}\|
          \sum_{j=1}^p \sqrt{V_{j,j}}
    \to 0 \qquad \text{as \ $n \to \infty$,}
 \end{align*}
 by \eqref{Gelfand2}, where for the last inequality we used \eqref{help8}.
This completes the proof.
\proofend

\section*{Acknowledgements}
We would like to thank the referees for their comments that helped us improve the paper.

\bibliographystyle{plain}
\bibliography{aggr5}

\def\polhk#1{\setbox0=\hbox{#1}{\ooalign{\hidewidth
  \lower1.5ex\hbox{`}\hidewidth\crcr\unhbox0}}}
\begin{thebibliography}{10}

\bibitem{BarIspPap2}
M.~Barczy, M.~Isp{\'a}ny, and G.~Pap.
\newblock Asymptotic behavior of unstable {${\rm INAR}(p)$} processes.
\newblock {\em Stochastic Process. Appl.}, 121(3):583--608, 2011.

\bibitem{BarNedPap}
M.~Barczy, F.~Ned\'enyi, and G.~Pap.
\newblock Iterated scaling limits for aggregation of randomized {INAR}(1)
  processes with idiosyncratic {P}oisson innovations.
\newblock {\em J. Math. Anal. Appl.}, 451(1):524--543, 2017.

\bibitem{DanPap}
Tivadar Danka and Gyula Pap.
\newblock Asymptotic behavior of critical indecomposable multi-type branching
  processes with immigration.
\newblock {\em ESAIM Probab. Stat.}, 20:238--260, 2016.

\bibitem{Gra}
C.~W.~J. Granger.
\newblock Long memory relationships and the aggregation of dynamic models.
\newblock {\em J. Econometrics}, 14(2):227--238, 1980.

\bibitem{HenSea}
Harold~V. Henderson and S.~R. Searle.
\newblock The vec-permutation matrix, the vec operator and {K}ronecker
  products: a review.
\newblock {\em Linear and Multilinear Algebra}, 9(4):271--288, 1980/81.

\bibitem{HornJohnson}
R.~A. Horn and Ch.~R. Johnson.
\newblock {\em Matrix analysis}.
\newblock Cambridge University Press, Cambridge, second edition, 2013.

\bibitem{JacShi}
J.~Jacod and A.~N. Shiryaev.
\newblock {\em Limit {T}heorems for {S}tochastic {P}rocesses}, volume 288 of
  {\em Grundlehren der Mathematischen Wissenschaften [Fundamental Principles of
  Mathematical Sciences]}.
\newblock Springer-Verlag, Berlin, second edition, 2003.

\bibitem{Jir}
M.~Jirak.
\newblock Limit theorems for aggregated linear processes.
\newblock {\em Adv.\ in Appl.\ Probab.}, 45(2):520--544, 2013.

\bibitem{KesSti2}
H.~Kesten and B.~P. Stigum.
\newblock Additional limit theorems for indecomposable multidimensional
  {G}alton-{W}atson processes.
\newblock {\em Ann. Math. Statist.}, 37:1463--1481, 1966.

\bibitem{KesSti3}
H.~Kesten and B.~P. Stigum.
\newblock Limit theorems for decomposable multi-dimensional {G}alton-{W}atson
  processes.
\newblock {\em J. Math. Anal. Appl.}, 17:309--338, 1967.

\bibitem{Latour}
Alain Latour.
\newblock Existence and stochastic structure of a non-negative integer-valued
  autoregressive process.
\newblock {\em J. Time Ser. Anal.}, 19(4):439--455, 1998.

\bibitem{Ned}
F.~Ned{\'e}nyi.
\newblock Conditional least squares estimators for multitype {G}alton--{W}atson
  processes.
\newblock {\em Acta Sci.\ Math.\ (Szeged)}, 81(1-2):325--348, 2015.

\bibitem{PilSur}
V.~Pilipauskait{\.e} and D.~Surgailis.
\newblock Joint temporal and contemporaneous aggregation of random-coefficient
  {AR}(1) processes.
\newblock {\em Stochastic Process.\ Appl.}, 124(2):1011--1035, 2014.

\bibitem{Quine}
M.~P. Quine.
\newblock The multi-type {G}alton-{W}atson process with immigration.
\newblock {\em J. Appl. Probability}, 7:411--422, 1970.

\bibitem{QuineDurham}
M.~P. Quine and P.~Durham.
\newblock Estimation for multitype branching processes.
\newblock {\em J. Appl. Probability}, 14(4):829--835, 1977.

\bibitem{Rob}
P.~M. Robinson.
\newblock Statistical inference for a random coefficient autoregressive model.
\newblock {\em Scand.\ J. Statist.}, 5(3):163--168, 1978.

\end{thebibliography}

\end{document}